\documentclass[aos,authoryear]{imsart}
\usepackage{amsmath,amssymb,amsthm,natbib,xspace,pstricks,pst-plot,xr}
\RequirePackage[colorlinks,citecolor=blue,urlcolor=blue]{hyperref}
\usepackage[all]{xy}
\setattribute{journal}{name}{}
\setattribute{infoline}{text}{}
\startlocaldefs
\newcommand{\IR}{\mathbb{R}}
\newcommand{\IN}{\mathbb{N}}
\newcommand{\ds}{\displaystyle}
\newcommand{\intersection}{\cap}
\newcommand{\Intersection}{\bigcap}
\newcommand{\union}{\cup}
\newcommand{\Union}{\bigcup}
\newcommand\PtwD{\textbf{\textup{(PtwD)}}\xspace}
\newcommand\ComD{\textbf{\textup{(ComD)}}\xspace}
\newcommand\UniD{\textbf{\textup{(UniD)}}\xspace}
\newcommand\PtwR{\textbf{\textup{(PtwR)}}\xspace}
\newcommand\ComR{\textbf{\textup{(ComR)}}\xspace}
\newcommand\RC{\textbf{\textup{(RC)}}\xspace}
\DeclareMathOperator{\D}{D}
\DeclareMathOperator{\ZD}{ZD}
\DeclareMathOperator{\HD}{HD}
\DeclareMathOperator{\MD}{MD}
\DeclareMathOperator{\AMD}{AMD}
\DeclareMathOperator{\WMD}{WMD}
\DeclareMathOperator{\interior}{int}
\DeclareMathOperator{\cl}{cl}
\DeclareMathOperator{\diam}{diam}
\DeclareMathOperator*{\Hlim}{H-lim}
\newtheorem{satz}{Theorem}[section]
\newtheorem{proposition}[satz]{Proposition}
\newtheorem{definition}{Definition}[section]
\newtheorem{corollary}{Corollary}[section]

\newtheorem{example}{Example}[section]
\allowdisplaybreaks

\newcommand{\manuallabel}[2]{\def\@currentlabel{#2}\label{#1}}
\endlocaldefs
\begin{document}
\begin{frontmatter}
\title{Convergence of depths and depth-trimmed regions}
\runtitle{Convergence of depths}
\begin{aug}
\author{\fnms{Rainer} \snm{Dyckerhoff}\ead[label=e1]{rainer.dyckerhoff@statistik.uni-koeln.de}
\ead[label=u1,url]{http://www.wisostat.uni-koeln.de/institut/professoren/dyckerhoff/}}
\runauthor{R. Dyckerhoff}
\affiliation{University of Cologne\thanksmark{m1}}
\address{Institute of Econometrics and Statistics\\
University of Cologne\\
Albertus-Magnus-Platz\\
50923 Cologne\\
Germany\\
\printead{e1}\\
\printead{u1}}
\end{aug}
\begin{abstract}
Depth is a concept that measures the `centrality' of a point in a given
data cloud or in a given probability distribution. Every depth defines a family 
of so-called trimmed regions. For statistical applications it is desirable that 
with increasing sample size the empirical depth as well as the empirical trimmed 
regions converge almost surely to their population counterparts.

In this article the connections between different types of convergence are
discussed. In particular, conditions are given under which the pointwise (resp. 
uniform) convergence of the data depth implies the pointwise (resp. compact) 
convergence of the trimmed regions in the Hausdorff metric as well as conditions 
under which the reverse implications hold. Further, it is shown that under 
relative weak conditions the pointwise convergence of the data depth (resp. 
trimmed regions) is equivalent to the uniform convergence of the data depth 
(resp. compact convergence of the trimmed regions). 
\end{abstract}


\end{frontmatter}
\section{Introduction}\label{sec.introduction}
In recent years data depth has been increasingly studied and is more and more
used in multivariate statistics. 
Applications of data depth in multivariate statistics include 
the construction of multivariate rank tests \citep{Liu92,LiuS93,Dyckerhoff02},
development of multivariate control charts \citep{Liu95},
construction of confidence regions \citep{YehS97}, 
multivariate data analysis \citep{LiuPS99},
cluster analysis \citep{Hoberg00,Hoberg03}, 
outlier detection \citep{Cramer03},
multivariate risk measurement \citep{CascosM07}, 
classification \citep{MoslerH06,LangeMM14},
and robust linear programming \citep{MoslerB14}.

Data depth is a function which quantifies the `centrality' of a point in a
given probability distribution. Closely related to the notion of depth is the
notion of \emph{central regions} or \emph{depth-trimmed regions}. Every depth 
defines a family of central regions in the following way. The $\alpha$-trimmed
region consists of all points that have a depth of at least $\alpha$ w.r.t. a 
given distribution. This is the set of points that have a certain degree of 
centrality and thus, this set is also called the $\alpha$-trimmed region. 
Since every depth defines a family of central regions and vice versa, the 
concepts of depth and central regions are in a sense equivalent.

Many depths have been proposed in the literature, e.g., 
the Mahalanobis depth \citep{Mahalanobis36}, 
halfspace depth \citep{Tukey75}, 
simplicial depth \citep{Liu88,Liu90}, 
majority depth \citep{Singh91}, 
projection depth \citetext{\citealp{Liu92,ZuoS00a,Zuo03}; based on a notion of 
outlyingness proposed by \citealp{Stahel81,Donoho82}},
zonoid depth \citep{KoshevoyM97},
weighted-mean depth \citep{DyckerhoffM11,DyckerhoffM12}
and others. These depths differ in many aspects, particularly
in the shape of trimmed regions or the deepest point. However, they share
certain properties which can be seen as desirable properties every depth should
satisfy. We define a depth as a function that satisfies certain postulates
which are stated in \cite{Dyckerhoff04}. Slightly differing 
sets of postulates have been given in \cite{Liu90} and \cite{ZuoS00a}.

In statistical applications the empirical depth, i.e., the depth w.r.t. the 
empirical measure defined by a sample $X_1,\dots,X_n$, is used as an estimator
for the depth w.r.t. the underlying distribution $P_X$. The same holds for the
empirical depth-trimmed regions. So the question of almost sure convergence of 
the empirical quantities to their population counterparts is of crucial 
interest since it is equivalent to strong consistency of these estimators. 

The convergence of depths and depth-trimmed regions has been extensively studied 
in the literature, e.g., for the the halfspace depth  
\citep{Eddy85,DonohoG92,Nolan92,MasseT94,Masse02,Masse04}, 
for the simplicial depth \citep{Liu90,Duembgen92}, 
for the zonoid depth \citep{KoshevoyM97,CascosLD16}, 
for the $\alpha$-trimming \citep{CascosLD08},
for the projection depth \citep{ZuoS00b,Zuo06}, 
for general type D depth functions \citep{ZuoS00b}, 
for generalized quantile functions defined by depth-trimmed regions \citep{Serfling02b}, 
for weighted-mean trimmed regions \citep{DyckerhoffM12}.

Results for general depths (but mainly for elliptically contoured distributions)
can be found in \cite{HeW97}. A generalization of these results for unimodal 
distributions with uniformly bounded and positive everywhere density is given by
\citep{Kim00}. For general depths and without assumptions on the
distributions first results on the connection between convergence of depths and
convergence of trimmed regions have been established in \cite{ZuoS00b}.

In the current paper we extend the results of \cite{ZuoS00b}. In particular,
we consider neither  
special depths nor special distributions. Instead, the results hold for all 
depths that satisfy the postulates of \cite{Dyckerhoff04}. Further we pose no 
restrictions on the considered distributions such as ellipticity or unimodality.
In particular we answer the following questions: Under what conditions does
pointwise (resp. uniform) convergence of the depth functions imply pointwise 
(resp. compact) convergence of the depth-trimmed regions and vice versa? 
Under what conditions does pointwise convergence of the depth functions (resp.
trimmed regions) imply uniform convergence of the depth functions (resp. compact
convergence of the trimmed regions)?

The paper is organized as follows. In Section~\ref{sec.concept} we
define data depth in an abstract way as a function that satisfies
a certain set of axioms. Section~\ref{sec.continuity} contains two results
on the continuity of the depth function and the trimmed regions. The main
results on convergence of depths and depth-trimmed regions as well as some
applications are then given in Section~\ref{sec.konv}. 
The proofs of the main results are collected in Appendix~\ref{sec.proofs}.
Some important results on Hausdorff convergence of sets are stated in 
Appendix~\ref{sec.Hausdorff}.  

In this paper we use the following notation. 
The complement of a set $A$ is denoted by $A^c$. 
Interior, closure, and boundary of a set $A$ are denoted by $\interior A$,
$\cl A$ and $\partial A$, respectively.
The Hausdorff distance of two non-empty compact sets $A$ and $B$ is denoted by 
$\delta_H(A,B)$. The Hausdorff-limit of a sequence $(A_n)_{n\in\IN}$
of non-empty compact sets is denoted by $\Hlim_{n\to\infty}A_n$.

\section{A general concept of data depth}\label{sec.concept}
We consider the depth of a point \emph{w.r.t. a probability distribution}.
Let ${\cal M}_0$ be the set of all probability measures on $(\IR^d,{\cal B}^d)$,
where ${\cal B}^d$ denotes the Borel $\sigma$-algebra on $\IR^d$,
and $\cal M$ a subset of ${\cal M}_0$. A depth assigns to each probability
measure $P\in{\cal M}$ a real function $\D(\cdot\,|\,P):\IR^d\to\IR_+$, the
so-called \emph{depth function} w.r.t. $P$. The set of all points that have a
depth of at least $\alpha$ is called the \emph{$\alpha$-trimmed region} or
\emph{$\alpha$-central region}. The $\alpha$-trimmed region w.r.t. $P$ is
denoted by $\D_\alpha(P)$, i.e., $\D_\alpha(P)=\{z\in\IR^d\,|\,\D(z\,|\,P)\ge\alpha\}$.

Often, the probability measure is the distribution $P_X$ of a $d$-variate random 
vector $X$. 
Since every probability measure $P$ on $(\IR^d,{\cal B}^d)$ can be represented 
as the distribution of a $d$-variate random vector $X$, every statement about 
depths can either be expressed in terms of probability measures or in terms of
random vectors. 

We now state some axioms that every reasonable notion of depth should satisfy.
\begin{description}
\item[D1: Affine invariance.] For every regular $d\times d$-Matrix $A$ and 
$b\in\IR^d$ holds $\D(z\,|\,P)=\D(Az+b\,|\,P_{Ax+b})$, where $P_{Ax+b}$ denotes
the image measure of $P$ under the transformation $x\mapsto Ax+b$.
\item[D2: Vanishing at infinity.] 
$\lim_{\|z\|\to\infty}\D(z\,|\,P)=0$.
\item[D3: Upper semicontinuity.] For each $\alpha>0$ the set $D_\alpha(P)$
is closed.
\item[D4: Monotone on rays.] For each $x_0$ of maximal depth and each
$r\in \IR^d$, $r\ne0$, the function 
$\lambda\mapsto \D(x_0+\lambda r\,|\,P)$, $\lambda\ge0$, is monotone decreasing.
\item[D4$^\prime$: Quasiconcavity.] For every $\alpha\ge 0$ the set 
$\D_\alpha(P)$ is convex.
\end{description}
The properties D1, D2, and D4 have been introduced by \cite{Liu90}. A further
set of axioms for a depth has been given by \cite{ZuoS00a}.
The main difference between their axioms and ours is that they do not require a
depth to be upper semicontinuous. In addition, they require that for
distributions having a properly defined unique center, the depth attains it
maximum value at this center. However, for centrally symmetric distributions,
this follows already from our axioms. 
For a discussion of these axioms, see e.g., \cite{Dyckerhoff04}.

\begin{definition}\label{def.datentiefe}
A mapping $\D$, that assigns to each probability measure $P$ in a certain set 
$\cal M$ a function $\D(\cdot\,|\,P):\IR^d\to\IR$ and that satisfies
the properties D1, D2, D3 and D4 is called \emph{depth}.
A depth that satisfies D4$^\prime$ is called \emph{convex depth}.
\end{definition}
A depth always attains its maximum on $\IR^d$. We denote this maximum depth
by $\alpha_{\max}(P)=\max\{\D(z\,|\,P)\,|\,z\in\IR^d\}$. 
A depth that has the same maximum depth for all probability measures $P$ is 
called a \emph{normed depth}. 

Properties D1 to D4 are formulated in terms of the depth itself. However, these
properties can also be formulated in terms of the trimmed regions. We now
state these equivalent properties. 
\begin{description}
\item[R1: Affine equivariance.] For every regular $d\times d$-matrix $A$ and
$b\in\IR^d$ holds $\D_\alpha(P_{Ax+b})=A\D_\alpha(P)+b$.
\item[R2: Boundedness.] For every $\alpha>0$ the region $\D_\alpha(P)$ is bounded.
\item[R3: Closedness.] For every $\alpha>0$ the region $\D_\alpha(P)$ is closed.
\item[R4: Starshapedness.] If $x_0$ is contained in all non-empty regions 
$\D_\alpha(P)$, then the non-empty regions $\D_\alpha(P)$ are 
starshaped w.r.t. $x_0$.
\item[R4$^\prime$: Convexity.] For every $\alpha>0$ the region $\D_\alpha(P)$ 
is convex.
\end{description}

In \cite{Dyckerhoff04} it has been shown that each of the statements 
\textup{\textbf{D1}, \textbf{D2}, \textbf{D3}, \textbf{D4}, \textbf{D4$^\prime$}}
is equivalent to the corresponding statement \textup{\textbf{R1}, \textbf{R2},
\textbf{R3}, \textbf{R4}, \textbf{R4$^\prime$}}.
A further important property that is satisfied by the trimmed regions of a depth
is left-continuity of the trimmed regions.
\begin{description}
\item[R5: Left-continuity.] For every $\alpha>0$ holds
$\D_\alpha(P)=\Intersection_{\beta:\beta <\alpha}\D_\beta(P)$.
Further, $\Intersection_{\alpha:\alpha\ge0}\D_\alpha(P)=\emptyset$ for
every $P\in{\cal M}$.
\end{description}
\textbf{R5} has been called \emph{intersection property} in \cite{Dyckerhoff04}.

In particular, it follows from R5 that the $\alpha$-trimmed
regions are monotone decreasing in $\alpha$, i.e., for $0\le\alpha_1\le\alpha_2$
holds $\D_{\alpha_2}(P)\subset \D_{\alpha_1}(P)$.

At the beginning of this section we have started from a depth and
defined the trimmed regions through the relation
$\D_\alpha(P):=\{z\,|\,\D(z\,|\,P)\ge\alpha\}$. However, one can also 
start from a family of trimmed regions and define the associated depth by its 
trimmed regions. 
In \cite{Dyckerhoff04} it has been shown that a family 
$(Z_\alpha(P))_{\alpha>0}$ of subsets of $\IR^d$ that satisfies the properties 
R1 to R5 defines a depth in the sense of Definition~\ref{def.datentiefe} via
\[                   
\D(z\,|\,P):=\sup\{\alpha\,|\,z\in Z_\alpha(P)\}\,.     
\]
If the sets $Z_\alpha(P)$ satisfy also R4$^\prime$, then $\D$ is a convex depth.
 
The idea of generating a depth by a suitable family of nested regions goes back
to \cite{Barnett76} and \cite{Eddy85}. It has already been used by
\cite{KoshevoyM97} to define the zonoid depth (see
Example~\ref{ex.ZD} below) and by \cite{Serfling02a} to define quantile
functions. By the above construction the weighted-mean depth can be constructed
from the weighted-mean trimmed regions defined in 
\cite{DyckerhoffM11,DyckerhoffM12}.

In the following we give some examples of depths which have already been proposed
in the literature.

\begin{example}\label{ex.MD}
The Mahalanobis depth \citep{Mahalanobis36} of a point $z$ is defined by
$\MD(z\,|\,P)=\left[1+(z-\mu_P)^\prime\Sigma_P^{-1}(z-\mu_P)\right]^{-1}$,
where $\mu_P$ denotes the expectation of $P$ and $\Sigma_P$ the covariance
matrix of $P$. It is well known that the Mahalanobis depth is a normed convex 
depth in the sense of Definition~\ref{def.datentiefe}. 
\end{example}

\begin{example}\label{ex.HD}
The halfspace depth \citep{Tukey75} is defined by
\[
\HD(z\,|\,P)=\inf\left\{P(H)\,|\,\text{$H$ is a closed halfspace containing $z$}
\right\}\,.
\]
The halfspace depth is a convex depth in the sense of 
Definition~\ref{def.datentiefe}.
A thorough discussion of the properties of the halfspace depth can be found 
in \cite{RousseeuwR99}.
\end{example}

\begin{example}\label{ex.ZD}
For a probability measure $P$ with finite first moments and $0<\alpha\le1$
the \emph{$\alpha$-zonoid trimmed region} $\ZD_\alpha(P)$ is defined by
\[
\ZD_\alpha(P)=\left\{\int\!\!xg(x)\,dP\,|\,g:\IR^d\to\left[0,\frac1\alpha\right]
\text{ measurable with} \int\!\!g(x)\,dP=1\right\}.
\]
The \emph{zonoid depth} is then defined by 
$\ZD(z\,|\,P)=\sup\{\alpha\,|\,z\in\ZD_\alpha(P)\}$.
The zonoid depth has been introduced by \cite{KoshevoyM97}.
It is a normed convex depth in the sense of Definition~\ref{def.datentiefe}.  
The properties of the zonoid depth and the associated zonoid trimmed regions are 
discussed in \cite{KoshevoyM97,Mosler02}.
\end{example}

\begin{example}\label{ex.WMD}
For a $d$-variate random vector $X$, the weighted-mean regions $\WMD_\alpha(P_X)$ 
\citep{DyckerhoffM11,DyckerhoffM12} are defined as the unique convex bodies 
whose support functions are given by  
\[
h(p)=\int_0^1Q_{p^\prime X}(t)\,dr_\alpha(t)\,,\quad p\in\IR^d\,,
\]
where $Q_{p^\prime X}$ denotes the quantile function of $p^\prime X$ and 
$r_\alpha$ is a suitable weighting function. 
The \emph{weighted-mean depth}, defined by 
$\WMD(z\,|\,P_X)=\linebreak\sup\{\alpha\,|\,z\in\WMD_\alpha(P_X)\}$, 
is a normed convex depth in the sense of Definition~\ref{def.datentiefe}.  
\end{example}

\begin{example}
The simplicial depth \citep{Liu88,Liu90} and the majority depth
\citep{Singh91} are no depths in the sense of Definition~\ref{def.datentiefe}.
The simplicial depth fails to satisfy D4 for discrete distributions (see the
counterexample in \citealp{ZuoS00a}). However, restricted to the class of
angular symmetric distributions the simplicial depth is a depth in the sense of
Definition~\ref{def.datentiefe}. The majority depth does not satisfy D2.
\end{example}

\section{Continuity of depths and trimmed regions}\label{sec.continuity}
In this section we consider the continuity properties of depths and trimmed 
regions. By \textbf{D3}, any depth in the sense of 
Definition~\ref{def.datentiefe} is upper semicontinuous. The following
theorem characterizes the depths that are even continuous. 

\begin{satz}\label{th.contdepth}
Let $\D$ be a depth. Then, the mapping $z\mapsto \D(z\,|\,P)$, $z\in\IR^d$,
is continuous if and only if for all $\beta>\alpha$ holds
\[
\D_\beta(P)\subset\interior\,\D_\alpha(P)\,.
\]
\end{satz}
The proofs of all theorems in Sections~\ref{sec.continuity} and \ref{sec.konv}
can be found in Appendix~\ref{sec.proofs}.

Next, we consider the question, under which conditions the trimmed regions
are continuous. Since because of \textbf{R2} and \textbf{R3} the trimmed regions 
$\D_\alpha(P)$ are compact sets, we use the usual notion of convergence for such 
sets, i.e., convergence in the \emph{Hausdorff-metric} or short  
\emph{Hausdorff-convergence}. The definition as well as some important facts on 
Hausdorff convergence are given in Appendix~\ref{sec.Hausdorff}.

Since the $\alpha$-trimmed region is the intersection of all $\beta$-trimmed
regions with $\beta<\alpha$ (\textbf{R5}), it is not surprising 
that for every depth in the sense of Definition~\ref{def.datentiefe} the mapping 
$\alpha\mapsto\D_\alpha(P)$ is left-continuous. Now, under which conditions is 
this mapping right-continuous, too? Intuitively, one has to demand, that the 
$\beta$-trimmed region is not much smaller than the $\alpha$-trimmed region, 
whenever $\beta$ is not much larger than $\alpha$. The following theorem shows 
that the trimmed regions are continuous in $\alpha$ if and only if for every 
$\alpha$ the $\alpha$-trimmed region is the closure of all points $z$ with 
$\D(z)>\alpha$. This condition can be seen as some kind of strict monotonicity 
of the depth. If $z$ is a point of depth $\alpha\in\bigl(0,\alpha_{\max}(P)\bigr)$ 
then each neighborhood of $z$ contains points of depth larger than $\alpha$. In 
other words, there is no neighborhood of $z$ on which the depth is constant. 
Therefore, if for a probability measure $P$ a depth satisfies this property,
we will say that the depth is strictly monotone for $P$.

\begin{definition}\label{def.strengmonoton}
Let $\D$ be a depth and $P$ a probability measure. 
If for each $\alpha\in\bigl(0,\alpha_{\max}(P)\bigr)$ we have
\[
\D_\alpha(P)=\cl{\{z\in\IR^d\,|\,\D(z\,|\,P)>\alpha\}}\,,
\]
then $\D$ is said to be strictly monotone for $P$.
\end{definition}

\begin{satz}\label{th.contalpha}
Let $\D$ be a depth. 
Then, the following assertions hold:
\begin{enumerate}
\item The mapping $\alpha\mapsto \D_\alpha(P)$, $0<\alpha\le\alpha_{\max}(P)$,
is left-continuous w.r.t. the Hausdorff metric, i.e., for every sequence
$(\alpha_n)_{n\in\IN}$ such that $\alpha_n<\alpha_0$ and
$\lim_{n\to\infty}\alpha_n=\alpha_0$ we have
$\Hlim_{n\to\infty}\D_{\alpha_n}(P)=\D_{\alpha_0}(P)$.
\item The mapping $\alpha\mapsto \D_\alpha(P)$, $0<\alpha\le\alpha_{\max}(P)$,
is continuous w.r.t. the Hausdorff metric if and only if $\D$ is strictly
monotone for $P$.
\end{enumerate}
\end{satz}

\section{Convergence of depths and trimmed regions}\label{sec.konv}
In this section we consider the following problem. Let $P_1,P_2,\dots$
be a sequence of probability measures on $(\IR^d,{\cal B}^d)$  and $P$ a further
probability measure on $(\IR^d,{\cal B}^d)$. If the corresponding depth functions
$\D(\,\cdot\,|\,P_n)$ converge to the depth function $\D(\,\cdot\,|\,P)$, what 
can be said about convergence of the corresponding trimmed regions? 
Of course the question can be posed also the other way: If the trimmed regions
$\D_\alpha(P_n)$ converge to $\D_\alpha(P)$, what can be said about convergence
of the depths? In this paper we make no assumptions on the sequence
$P_1,P_2,\dots$ of probability measures. However, the main application of our
results is the following situation: 
Let $X$ be $d$-variate random vector with distribution $P_X$ and $X_1,X_2,\dots$ 
a sequence of random vectors that are independent and identically distributed 
with distribution $P_X$. Denote by $P_n$ the empirical distribution on $X_1,
\dots,X_n$, i.e., the distribution that assigns probability $1/n$ to each of 
these $n$ points. The \emph{empirical depth} $\D(\,\cdot\,|\,P_n)$ then 
constitutes an estimator for the population depth $\D(\,\cdot\,|\,P_X)$. 
In the same way the \emph{empirical trimmed regions} $\D_\alpha(P_n)$ are 
estimators for the population trimmed region $\D_\alpha(P_X)$. These estimators 
are strongly consistent if they converge with probability one to their 
population counterparts. We will comment on this situation later. 

Intuitively, one would assume that if the empirical depths converge, the
empirical trimmed regions will converge, too, and vice versa. We will show
that in most circumstances this is true. However, there are also situations
where the depths converge pointwise or even uniform but the corresponding
trimmed regions do not converge for certain values of $\alpha$. Conversely,
it may also happen that the trimmed regions do converge whereas the depths
themselves do not. If one thinks deeper about this question, this is not
surprising. Consider for example the sequence of empirical distribution
functions. This sequence converges with probability one uniformly to the
theoretical distribution function, whereas the sequence of empirical quantile
functions need in general not converge for all $p$.
The exact relationships between convergence of depths and convergence of
corresponding trimmed regions are in fact far more complicated than this is the
case for the relationship between convergence of empirical distribution and
quantile functions.

We investigate the connections between the following notions of convergence
for depths and depth-trimmed regions.
\begin{description}
\item[\PtwD Pointwise convergence of depths]~\\ For all $z\in\IR^d$ holds
$\ds\lim_{n\to\infty}\D(z\,|\,P_n)=\D(z\,|\,P)$.
\item[\ComD Compact convergence of depths]~\\ For every compact set
$M\subset\IR^d$ holds
\[
\lim_{n\to\infty}\sup_{z\in M}|\D(z\,|\,P_n)-\D(z\,|\,P)|=0\,.
\]
\item[\UniD Uniform convergence of depths]~\\ It holds
\[
\lim_{n\to\infty}\sup_{z\in \IR^d}|\D(z\,|\,P_n)-\D(z\,|\,P)|=0\,.
\]
\item[\PtwR Pointwise Hausdorff-convergence of trimmed regions]~\\
For every $\alpha\in\bigl(0,\alpha_{\max}(P)\bigr)$ holds 
$\ds\Hlim_{n\to\infty}\D_\alpha(P_n)=\D_\alpha(P)$.
\item[\ComR Compact Hausdorff-convergence of trimmed regions]~\\
For every compact interval $A\subset\bigl(0,\alpha_{\max}(P)\bigr)$ holds
\[
\lim_{n\to\infty}\sup_{\alpha\in A}\delta_H(\D_\alpha(P_n),\D_\alpha(P))=0\,.
\]
\end{description}
A further complication results from the following observation. The Hausdorff
distance between two sets is undefined whenever one of the sets is empty. For
this reason we have considered only values of $\alpha$ in 
$\bigl(0,\alpha_{\max}(P)\bigr)$
in \PtwR and \ComR. For these values of $\alpha$ at least $\D_\alpha(P)$ is
non-empty. Of course it cannot be guaranteed that for these $\alpha$'s the
trimmed regions $\D_\alpha(P_n)$ are non-empty for each $n$, too. However,
if we assume \PtwR and \ComR, it is implicitly assumed that for sufficiently
large $n$ the distances $\delta_H(\D_\alpha(P_n),\D_\alpha(P))$ and
$\sup_{\alpha\in A}\delta_H(\D_\alpha(P_n),\D_\alpha(P))$ are defined.
For finitely many $n$ these distances may be undefined.

Conversely it may occur that for every $n$ the maximum depth $\alpha_{\max}(P_n)$ 
w.r.t. $P_n$ is greater than the maximum depth $\alpha_{\max}(P)$ w.r.t. $P$. 
In this case \PtwR and \ComR make assertions on the convergence of trimmed regions
$\D_\alpha(P_n)$ for $\alpha\in\bigl(0,\alpha_{\max}(P)\bigr)$, but not on the 
`convergence' of the trimmed regions $\D_\alpha(P_n)$ for 
$\alpha>\alpha_{\max}(P)$. 
Therefore, we occasionally need a further condition that guarantees that the 
trimmed regions $\D_\alpha(P_n)$ converge for $\alpha>\alpha_{\max}(P)$ to the 
empty set:
\begin{description}
\item[\RC Range condition] 
It holds $\ds\limsup_{n\to\infty}\alpha_{\max}(P_n)\le\alpha_{\max}(P)$.
\end{description}
We will use this condition mostly in the following equivalent form:
\begin{description}
\item[\RC Range condition] 
For each $\alpha>\alpha_{\max}(P)$ there exists $N_\alpha\in \IN$, such that
$\D_\alpha(P_n)=\emptyset$ for all $n\ge N_\alpha$.
\end{description}
The condition \RC is trivially satisfied if the maximum depth is the same
for all distributions, i.e., for normed depths. 
Examples of normed depths are the Mahalanobis depth and the zonoid depth. 
For both of these depths the point of maximum depth has always depth one.

\textbf{General assumption:} To avoid technical difficulties we assume in this 
section that the $\alpha$-trimmed regions have full dimension for 
$0<\alpha<\alpha_{\max}(P)$. This condition is satisfied by the commonly used 
depths, unless the probability measure $P$ is concentrated on a hyperplane.

We first prove a theorem that relates the pointwise convergence of depths
to the set-theoretic limit of the trimmed regions. 
The limit inferior and the limit superior of a sequence of sets are defined by
\[
\liminf_{n\to\infty}A_n=\Union_{n=1}^\infty\Intersection_{m=n}^\infty A_m
\quad\text{resp.}\quad
\limsup_{n\to\infty}A_n=\Intersection_{n=1}^\infty\Union_{m=n}^\infty A_m\,.
\]
Clearly, the limit inferior is the set of all points that are eventually 
contained in all sets $A_n$, the limit superior is the set of all points that 
are contained in infinitely many of the sets $A_n$. 
Hence, $\liminf_{n\to\infty}A_n\subset\limsup_{n\to\infty}A_n$. 
If $\liminf_{n\to\infty}A_n=\limsup_{n\to\infty}A_n$, then the sequence 
$(A_n)_{n\in\IN}$ is said to converge (in the set-theoretic sense) and one 
defines
\[
\lim_{n\to\infty}A_n:=\liminf_{n\to\infty}A_n=\limsup_{n\to\infty}A_n\,.
\]

\begin{satz}\label{th.konvpunkt}
The following two statements are equivalent:
\begin{enumerate}
\item \PtwD
\item For every $\alpha\ge 0$ holds
\[
\{z\in\IR^d\,|\,\D(z\,|\,P)>\alpha\}\subset\liminf_{n\to\infty}\D_\alpha(P_n)
\subset\limsup_{n\to\infty}\D_\alpha(P_n)\subset\D_\alpha(P)\,.
\]
\end{enumerate}
\end{satz}

\textbf{Remark:} Since the set $\{z\in\IR^d\,|\,\D(z\,|\,P)>\alpha\}$ is
non-empty for $\alpha\in\bigl(0,\alpha_{\max}(P)\bigr)$ the set 
$\liminf_{n\to\infty}\D_\alpha(P_n)$ is non-empty as well. In particular, it 
follows from the preceding theorem that the trimmed regions $\D_\alpha(P_n)$ are 
eventually non-empty for every $\alpha\in\bigl(0,\alpha_{\max}(P)\bigr)$.

The preceding theorem shows that in general the set-theoretic convergence of
trimmed regions does not follow from the pointwise convergence of depths.
This conclusion is justified only when $\{z\in\IR^d\,|\,\D(z\,|\,P)>\alpha\}=
\D_\alpha(P)$. However, if the depth is continuous for $P$, this condition is
never satisfied.

If one considers not the set-theoretic convergence but the Hausdorff convergence 
of trimmed regions the following theorem can be deduced.
It says that for convex and strictly monotone depths the pointwise convergence
of depths implies the Hausdorff convergence of the trimmed regions.
\begin{satz}\label{th.PktTPktZ}
Let $\D$ be a convex depth and let $\D$ be strictly monotone for $P$.
Then, \PtwD $\implies$ \PtwR.

If $\D_{\alpha_{\max}(P)}(P)$ is a singleton and if the trimmed regions 
$\D_{\alpha_{\max}(P)}(P_n)$ are eventually non-empty, then \PtwR holds also 
for $\alpha=\alpha_{\max}(P)$.
\end{satz}
As is shown in Example~\ref{ex.1} in the supplement \citep{Dyckerhoff17},
without the assumption of strict monotonicity the Hausdorff convergence of the
trimmed regions can in general not be concluded.

The following theorem shows that for continuous convex depths the pointwise
convergence of depths follows from the Hausdorff convergence of the trimmed 
regions.

\begin{satz}\label{th.PktZPktT}
Let $\D$ be a convex depth and let $\D$ be continuous for $P$.
Then, \emph{[}\PtwR and \RC\emph{]} $\implies$ \PtwD.
\end{satz}
An example that shows that without the assumption of continuity the above
theorem is in general false is given in Example~\ref{ex.3} in the supplement
\citep{Dyckerhoff17}.

We will now study what are the implications of compact or uniform convergence of
depths. On the one hand the following theorem shows that compact and uniform
convergence of depths are in fact equivalent. On the other hand it gives two
conditions on the trimmed regions that are equivalent to the uniform convergence
of depths. Essentially, these two conditions state that for sufficiently large
$n$ the empirical trimmed regions $\D_\alpha(P_n)$ lie between the trimmed regions
$\D_{\alpha+\epsilon}(P)$ and $\D_{\alpha-\epsilon}(P)$. Further, it will be shown
that -- in contrast to the pointwise convergence -- the uniform convergence
implies the condition \RC.

\begin{satz}\label{th.konvglm}
The following statements are equivalent:
\begin{enumerate}
\item \UniD
\item \ComD
\item For every $\alpha> 0$ and every $\epsilon>0$ there exists an
$N_{\alpha,\epsilon}\in\IN$, such that
\[
\D_{\beta+\epsilon}(P)\subset\D_\beta(P_n)\subset\D_{\beta-\epsilon}(P)
\quad\text{for all $n\ge N_{\alpha,\epsilon}$ and $\beta\ge\alpha$.}
\]
\item It holds \RC and for every compact interval 
$A\subset\bigl(0,\alpha_{\max}(P)\bigr)$ and $\epsilon>0$ there exists an 
$N_{A,\epsilon}\in\IN$, such that
\[
\D_{\beta+\epsilon}(P)\subset\D_\beta(P_n)\subset\D_{\beta-\epsilon}(P)
\quad\text{for all $n\ge N_{A,\epsilon}$ and $\beta\in A$.}
\]
\end{enumerate}
\end{satz}

The implication $(ii)\implies(iii)$ in Theorem~\ref{th.konvglm} was already
proved in Theorem~4.1 in \cite{ZuoS00b}. 

The following theorem shows that for strictly monotone depths the compact
convergence of the depth implies compact Hausdorff convergence of trimmed
regions. As in Theorem~\ref{th.PktTPktZ}, without the assumption of strict 
monotonicity this statement is in general not valid. This can be seen from
Example~\ref{ex.1} in the supplement \citep{Dyckerhoff17}. 

\begin{satz}\label{th.KomTKomZ}
Let $\D$ be strictly monotone for $P$. Then, \ComD $\implies$ \ComR.

If $\D_{\alpha_{\max}(P)}(P)$ is a singleton and if the trimmed regions 
$\D_{\alpha_{\max}(P)}(P_n)$ are eventually non-empty, then \ComR holds also for 
every compact interval $A\subset\bigl(0,\alpha_{\max}(P)\bigr]$.
\end{satz}

If the depth is convex and continuous for $P$ and if the condition \RC holds,
the converse of the preceding theorem holds, i.e., compact Hausdorff convergence
of the trimmed regions implies uniform convergence of the depths.
Again, as in Theorem~\ref{th.PktZPktT}, without the assumption of continuity 
this does in general not hold, see Example~\ref{ex.3} in the supplement 
\cite{Dyckerhoff17}. 
\begin{satz}\label{th.KomZGlmT}
Let $\D$ be a convex depth that is continuous for $P$. Then,
\emph{[}\ComR and \RC\emph{]} $\implies$ \UniD.
\end{satz}

We have seen in Theorem~\ref{th.konvglm} that for depths compact and uniform
convergence are equivalent. We will see in the next two theorems that
under mild conditions even pointwise and compact convergence are
equivalent. This holds for the convergence of depths as well as for the
Hausdorff convergence of trimmed regions. 
\begin{satz}\label{th.PktZKomZ}
Let $\D$ be strictly monotone for $P$. Then, \PtwR $\iff$ \ComR.

If $\D_{\alpha_{\max}(P)}(P)$ is a singleton, then the open interval 
$\bigl(0,\alpha_{\max}(P)\bigr)$ can be replaced by the half-open interval 
$\bigl(0,\alpha_{\max}(P)\bigr]$ in \PtwR and \ComR.
\end{satz}
The condition of strict monotonicity is crucial in the above theorem. Without 
strict monotonicity, Theorem~\ref{th.PktZKomZ} is not valid as can be seen
from Example~\ref{ex.2} in the supplement \citep{Dyckerhoff17}. 

We now state the result for convergence of depths. Here the depth has to be
convex and continuous for $P$. Again, as in Theorems~\ref{th.PktZPktT} 
and \ref{th.KomZGlmT} the additional condition \RC is needed to ensure that
pointwise and uniform convergence are equivalent.
\begin{satz}\label{th.PktTGlmT}
Let $\D$ be a convex depth that is continuous for $P$. Then,
\emph{[}\PtwD and \RC\emph{]}$\iff$\UniD.
\end{satz}
Example~\ref{ex.4} in the supplement \citep{Dyckerhoff17} shows that without 
the assumption of continuity the above theorem is in general false. 

The connections between the different notions of convergence can be illustrated
nicely by a diagram. The following figure shows the implications which result
from the preceding theorems as well as the corresponding assumptions. For better
clarity we have replaced the implication arrows by simple arrows.
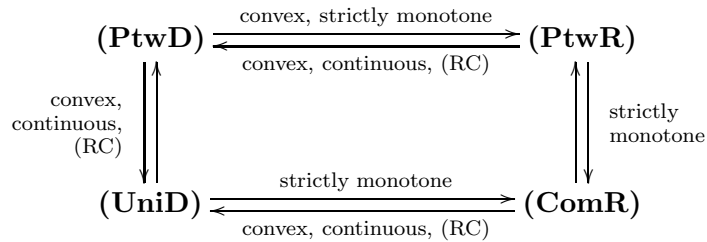
\[
{\def\ss{\scriptstyle}
\xymatrix @C=4cm @R=1.6cm {
\PtwD \ar@<0.5ex>[r]^{\text{convex, strictly monotone}}
             \ar@<-0.5ex>[d]_{\begin{array}{r}\ss\text{convex,}\\[-1ex]\ss\text{continuous,}\\[-1ex]\ss\text{(RC)}\end{array}} &
\PtwR \ar@<0.5ex>[l]^{\text{convex, continuous, (RC)}}
             \ar@<0.5ex>[d]^{\begin{array}{l}\ss\text{strictly}\\[-1ex]\ss\text{monotone}\end{array}}\\
\UniD \ar@<0.5ex>[r]^{\text{strictly monotone}}
             \ar@<-0.5ex>[u] &
\ComR \ar@<0.5ex>[l]^{\text{convex, continuous, (RC)}}
             \ar@<0.5ex>[u]
}}
\]
\textbf{Remark:} The condition that the depth be strictly monotone for $P$ can also
be seen as a continuity condition. In fact it follows from
Theorem~\ref{th.contalpha} that the mapping $\alpha\mapsto\D_\alpha(P)$ is
continuous in this case. Thus, one could as well replace the condition
`$\D$ is strictly monotone for $P$' by `the mapping $\alpha\mapsto\D_\alpha(P)$
is continuous'.

For the important class of normed and convex depths one gets the following 
connections:
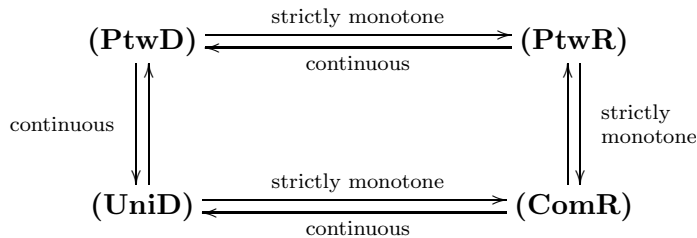
\[
{\def\ss{\scriptstyle}
\xymatrix @C=4cm @R=1.6cm {
\PtwD \ar@<0.5ex>[r]^{\text{strictly monotone}}
             \ar@<-0.5ex>[d]_{\begin{array}{r}\ss\text{continuous}\end{array}} &
\PtwR \ar@<0.5ex>[l]^{\text{continuous}}
             \ar@<0.5ex>[d]^{\begin{array}{l}\ss\text{strictly}\\[-1ex]\ss\text{monotone}\end{array}}\\
\UniD \ar@<0.5ex>[r]^{\text{strictly monotone}}
             \ar@<-0.5ex>[u] &
\ComR \ar@<0.5ex>[l]^{\text{continuous}}
             \ar@<0.5ex>[u]
}}
\]

A typical application of the above theorems arises when $X_1,X_2,\dots$ is a 
sequence of $d$-variate random vectors, defined on a joint probability space
$(\Omega,{\cal A},P)$, that are independent and identically distributed with 
distribution $P_X$, in symbols $X_1,X_2,\dots\stackrel{iid}{\sim}P_X$. 
Then, let $P_n$ be the empirical measure on $X_1,\dots,X_n$, i.e.,
$P_n=1/n\sum_{i=1}^n\varepsilon_{X_i}$, where $\varepsilon_{X_i}$ denotes the 
one-point measure on $X_i$. Note  that $P_n$ is in fact a random measure since 
it depends on the concrete realizations $X_1(\omega),\dots,X_n(\omega)$.
It is well known that with probability one the empirical measures converges
weakly to the distribution $P_X$. 
In this situation all of the above theorems have corollaries like the
following. 
\begin{corollary}[to Theorem~\ref{th.PktTGlmT}]\label{cor.PktTGlmT}
Let $X_1,X_2,\dots\stackrel{iid}{\sim}P_X$ and $P_n$ be the empirical measure
on $X_1,\dots,X_n$. Let further $\D$ be a convex depth that is continuous for 
$P_X$. Then,
\[
\text{\emph{[}\PtwD and \RC\emph{]} $P$-almost surely}
\iff
\text{\UniD $P$-almost surely.}
\]
\end{corollary}
Analogous corollaries hold for all of the above theorems. We do not state them 
here to avoid unnecessary repetitions.

Here, one has to be careful to distinguish between the probability measure
$P$ on the underlying probability space $(\Omega,{\cal A},P)$ and the 
probability measure $P_X$ that is the distribution of each of the random 
variables $X_i$. The depth is computed w.r.t. the distribution $P_X$, whereas
`$P$-almost surely' refers to the measure $P$ of the underlying probability
space.

We illustrate the application of the above results with some examples.
\begin{example}[Mahalanobis depth, see Example~\ref{ex.MD}]
The Mahalanobis depth is a normed convex depth. It is continuous and 
strictly monotone for each $P$.
From the strong law of large numbers follows (with probability one) the 
pointwise convergence of the empirical Mahalanobis depth.
Thus, with probability one, the empirical Mahalanobis depth converges
uniformly to its population version and the empirical $\alpha$-trimmed regions 
converge compactly on $(0,1)$.
\end{example}

\begin{example}[Halfspace depth, see Example~\ref{ex.HD}]
The halfspace depth is a convex depth. It is continuous for distributions with 
density. Under some additional assumptions on $P$ (e.g., convex support) it is 
also strictly monotone.
It is easy to show that with probability one the halfspace depth converges 
pointwise and the range condition is satisfied.
Thus, under the above conditions, with probability one, the empirical halfspace 
depth converges uniformly to its population version and the empirical 
$\alpha$-trimmed regions converge compactly on $\bigl(0,\alpha_{\max}(P)\bigr)$. 
\end{example}

\begin{example}[Zonoid depth, see Example~\ref{ex.ZD}]
The zonoid depth is a normed convex depth that is strictly monotone and
continuous for distributions with density.
It was shown in \cite{Mosler02} that, with probability one, the empirical zonoid 
regions converge pointwise to their population version.
Thus, for distributions with density, the empirical zonoid depth converges 
uniformly to its population version and the empirical $\alpha$-trimmed regions 
converge compactly on $(0,1)$.
\end{example}

\begin{example}[Asymmetric Mahalanobis depth]
The asymmetric Mahalanobis depth \cite[see][]{Dyckerhoff04} is defined by
\[
\AMD(z\,|\,P)=\inf_{p\in S^{d-1}}
\left[1+\left(
\frac{p^\prime z-\mu_{p^\prime P}}
     {\sigma^+_{p^\prime P}}
\right)\right]^{-1}\,,
\]
where ${\sigma^+_P}^2$ denotes the upper semi-variance of $P$ and $S^{d-1}$
is the unit sphere in $\IR^d$.
The asymmetric Mahalanobis depth is a normed convex depth. 
It is continuous and strictly monotone for each $P$.
From the strong law of large numbers follows (with probability one) the 
pointwise convergence of the empirical asymmetric Mahalanobis depth.
Thus, with probability one, the empirical asymmetric Mahalanobis depth converges
uniformly to its population version and the empirical $\alpha$-trimmed regions
converge compactly on $(0,1)$.
\end{example}

\begin{example}[Weighted-mean depth, see Example~\ref{ex.MD}]
The weighted-mean depth is a normed convex depth that is strictly monotone 
and continuous for distributions with density.
Let $(X_n)$ be a sequence of random vectors with finite first moments that
converges in distribution to a random vector $X$. It was shown in 
\cite{DyckerhoffM12} that the weighted-mean regions $\WMD_\alpha(X_n)$ are
pointwise Hausdorff convergent to $\WMD_\alpha(X)$ provided the sequence 
$(X_n)$ is uniformly integrable. Thus, it follows from the above theorems that 
the weighted-mean regions are even compact convergent on $(0,1]$. For  
distributions with density the associated depth functions converge uniformly
on $\IR^d$.
\end{example}

\begin{appendix}
\section{Proofs of the main theorems}\label{sec.proofs}
In this section we use the following notation: For a given depth $\D$ and
$n\in\IN$ we denote the $\alpha$-trimmed region $\D_\alpha(P_n)$ w.r.t. $P_n$
shortly by $\D_\alpha^n$ and the depth $\D(z\,|\,P_n)$ of a point $z\in\IR^d$
w.r.t. $P_n$ shortly with $\D^n(z)$. We use the notation $\D_\alpha$ for
$\D_\alpha(P)$ and $\D(z)$ for $\D(z\,|\,P)$ in the same way. In the same spirit 
we often write simply $\alpha_{\max}$ instead of $\alpha_{\max}(P)$. 

\textbf{Proof (of Theorem~\ref{th.contdepth}):} Assume that $\D_\beta\not\subset
\interior\,\D_\alpha$ for some $\beta>\alpha$. Then, there exists $z\in
\D_\beta$ that is contained in $\partial\D_\alpha$. Obviously,
$\D_\alpha\ne\IR^d$ in this case. Since $\D_\alpha$ is closed, its
complement is open and there exists a sequence $(z_n)_{n\in\IN}$ in
$\D_\alpha^c$ converging to $z$. But $\limsup_{n\to\infty}\D(z_n)\le
\alpha<\beta\le\D(z)$. Thus, the mapping $z\mapsto\D(z)$ is not
continuous.

Now assume that $\D_\beta\subset\interior\,\D_\alpha$ for all $\beta>
\alpha$. We show that $\D(\,\cdot\,|\,P)$ is lower semicontinuous, i.e., every
set $\{z\,|\,\D(z)>\alpha\}$ is open. For each 
$z_0\in\{z\,|\,\D(z)>\alpha\}$ there is $\gamma$ such that 
$\D(z_0)=:\beta>\gamma>\alpha$. Since $\D_\beta\subset\interior
\D_\gamma$, it follows $z_0\in\interior\D_\gamma$. Thus, there is a 
neighborhood $U$ of $z_0$ such that $U\subset\D_\gamma\subset
\{z\,|\,\D(z)>\alpha\}$ which shows that $\{z\,|\,\D(z)>\alpha\}$ 
is open. Hence, $\D(\,\cdot\,|\,P)$ is lower semicontinuous. Since 
$\D(\,\cdot\,|\,P)$ is also upper semicontinuous, $\D(\,\cdot\,|\,P)$ is 
continuous.\qed

\textbf{Proof (of Theorem~\ref{th.contalpha}):} We start with Part~$(i)$.
We assume w.l.o.g. that the sequence
$(\alpha_n)_{n\in\IN}$ is increasing. The sequence of sets $\D_{\alpha_n}$
is decreasing and it follows from Proposition~\ref{pr.HD-limes-aufab} in 
Appendix~\ref{sec.Hausdorff} that
\[
\Hlim_{n\to\infty}\D_{\alpha_n}=\Intersection_{n=1}^\infty\D_{\alpha_n}
=\D_{\alpha_0}\,.
\]
To prove Part~$(ii)$ it suffices to show that the strict monotonicity for $P$ is
equivalent to the right-continuity of the mapping $\alpha\mapsto \D_\alpha$.
This mapping is right continuous if and only if $\Hlim_{n\to\infty}
\D_{\alpha_n}=\D_{\alpha_0}$ for every sequence $(\alpha_n)_{n\in\IN}$ that
is decreasing to $\alpha_0\in\bigl(0,\alpha_{\max}(P)\bigr)$. Since the sequence 
of sets $\D_{\alpha_n}$ is increasing it follows from 
Proposition~\ref{pr.HD-limes-aufab} in Appendix~\ref{sec.Hausdorff} that
\[
\Hlim_{n\to\infty}\D_{\alpha_n}
=\cl\left(\Union_{n=1}^\infty\D_{\alpha_n}\right)
=\cl\left(\{z\in\IR^d\,|\,\D(z)>\alpha_0\}\right)\,.
\]
Thus, the mapping is right-continuous if and only if
\[
\cl\left(\{z\in\IR^d\,|\,\D(z)>\alpha_0\}\right)=\D_{\alpha_0}
\]
for each $\alpha_0\in\bigl(0,\alpha_{\max}(P)\bigr)$, i.e., if the depth is 
strictly monotone for $P$.
\qed

\textbf{Proof (of Theorem~\ref{th.konvpunkt}):}
We start with $(i)\implies(ii)$. Let $z\in D_{\alpha+\epsilon}$, then
$\D(z)\ge\alpha+\epsilon$. Since $\lim_{n\to\infty}\D^n(z)=\D(z)$ for every
$\epsilon>0$, there exists an $N_\epsilon\in\IN$ such that
\[
|\D^n(z)-\D(z)|<\epsilon\quad\text{for all $n\ge N_\epsilon$.}
\]
This implies
\[
\D^n(z)>\D(z)-\epsilon\ge\alpha+\epsilon-\epsilon=\alpha\,.
\]
Thus, $z\in\D_\alpha^n$ for all $n\ge N_\epsilon$ and therefore
$z\in\liminf_{n\to\infty}\D_\alpha^n$. This shows that
\begin{equation}\label{eq.konv1}
\D_{\alpha+\epsilon}\subset\liminf_{n\to\infty}\D_\alpha^n\quad\text{for all $\epsilon>0$.}
\end{equation}
Now we assume that $z\notin\D_{\alpha-\epsilon}$, i.e., $\D(z)<\alpha-\epsilon$.
From $(i)$ follows that there exists $N_\epsilon\in\IN$, such that
\[
|\D^n(z)-\D(z)|<\epsilon\quad\text{for all $n\ge N_\epsilon$.}
\]
Therefore,
\[
\D^n(z)<\D(z)+\epsilon<\alpha-\epsilon+\epsilon=\alpha\,.
\]
Thus, $z\in(\D_\alpha^n)^c$ for all $n\ge N_\epsilon$, and therefore
$z\in\liminf_{n\to\infty}(\D_\alpha^n)^c=(\limsup_{n\to\infty}\D_\alpha^n)^c$.
From this follows
\begin{equation}\label{eq.konv2}
\limsup_{n\to\infty}\D_\alpha^n\subset\D_{\alpha-\epsilon}\quad\text{for all $\epsilon>0$.}
\end{equation}
From the equations \eqref{eq.konv1} and \eqref{eq.konv2} follows
\[
\{z\,|\,\D(z)>\alpha\}=\Union_{\epsilon>0}\D_{\alpha+\epsilon}
\subset\liminf_{n\to\infty}\D_\alpha^n\subset\limsup_{n\to\infty}\D_\alpha^n\subset
\Intersection_{\epsilon>0}\D_{\alpha-\epsilon}=\D_\alpha\,,
\]
as was to be shown.

We now prove the direction $(ii)\implies(i)$. Let $z\in\IR^d$ such that
$\D(z)=\alpha$. We assume that the sequence $(\D^n(z))_{n\in\IN}$ does not
converge to $\alpha$. Then there is an $\epsilon>0$, such that
$|\D^n(z)-\D(z)|\ge\epsilon$ infinitely often. Thus we have
$\D^n(z)\ge\alpha+\epsilon$ or $\D^n(z)\le\alpha-\epsilon$ for infinitely many $n$.
In the first case $z\in\limsup_{n\to\infty}\D_{\alpha+\epsilon}$. From $(ii)$
follows $z\in\D_{\alpha+\epsilon}$, i.e., $\D(z)\ge \alpha+\epsilon$ in
contradiction to $\D(z)=\alpha$. In the second case
$z\in\limsup_{n\to\infty}(\D_{\alpha-\epsilon/2})^c=
(\liminf_{n\to\infty}\D_{\alpha-\epsilon/2})^c$. From $(ii)$ follows
$z\notin \{x\in\IR^d\,|\,\D(x)>\alpha-\epsilon/2\}$, i.e.,
$\D(z)\le\alpha-\epsilon/2$, in contradiction to $\D(z)=\alpha$.\qed

\textbf{Proof (of Theorem~\ref{th.PktTPktZ}):} We show that \PtwD implies that for every $\alpha
\in(0,\alpha_{\max})$ and for every $M>0$ the equation
\begin{equation}\label{eq.konv4}
\lim_{n\to\infty}\max\{\delta(x,\D_\alpha^n)\,|\,x\in\D_\alpha\intersection B(0,M)\}=0
\end{equation}
as well as the equation
\begin{equation}\label{eq.konv5}
\lim_{n\to\infty}\max\{\delta(x,\D_\alpha)\,|\,x\in\D_\alpha^n\intersection B(0,M)\}=0
\end{equation}
hold. Since the trimmed regions $\D_\alpha^n$ are connected, it then follows
from Theorem~\ref{th.AWHD} that $\Hlim_{n\to\infty}\D_\alpha^n=\D_\alpha$.

To show \eqref{eq.konv4} we first show the slightly stronger assertion
\begin{equation}\label{eq.konv6}
\lim_{n\to\infty}\max\{\delta(x,\D_\alpha^n)\,|\,x\in\D_\alpha\}=0\,.
\end{equation}
Obviously, \eqref{eq.konv6} implies \eqref{eq.konv4}.   
If \eqref{eq.konv6} does not hold, then there
exists $\epsilon>0$ and a subsequence $(x_{n_k})_{k\in\IN}$
with $x_{n_k}\in \D_\alpha$ such that $\delta(x_{n_k},\D_\alpha^{n_k})>\epsilon$
for all $k\in\IN$. Since $\D_\alpha$ is compact, the sequence $(x_{n_k})_{k\in\IN}$
has a convergent subsequence. We therefore assume w.l.o.g. that the sequence
$(x_{n_k})$ itself is convergent with $\lim_{k\to\infty}x_{n_k}=x_0\in\D_\alpha$.
For sufficiently large $k$ we have $\|x_0-x_{n_k}\|<\frac\epsilon2$ and
\[
\delta(x_0,\D_\alpha^{n_k})\ge \delta(x_{n_k},\D_\alpha^{n_k}) - \|x_0-x_{n_k}\|
> \epsilon - \frac\epsilon2 = \frac\epsilon2\,,
\]
i.e.,
\[
B\Bigl(x_0,\frac\epsilon2\Bigr)\intersection \D_\alpha^{n_k}=\emptyset\,.
\]
Since $\D_\alpha$ is the closure of all points with depth greater than $\alpha$,
there exists a point $z\in B(x_0,\frac\epsilon2)$ with $\D(z)>\alpha$.
It follows from Theorem~\ref{th.konvpunkt} that $z\in\liminf_{n\to\infty}
\D_\alpha^n$, i.e., the sets $\D_\alpha^n$ eventually contain $z$.
On the other hand $z\notin\D_\alpha^{n_k}$ for infinitely many $k$,
contradiction. Thus, \eqref{eq.konv6} and therefore also \eqref{eq.konv4} holds.

We now show that \eqref{eq.konv5} holds for every $M>0$. Assume that this is not 
the case.  
Then there is an $M>0$ and an $\epsilon>0$ as well as a sequence
$(x_{n_k})_{k\in\IN}$ with $\|x_{n_k}\|\le M$ such that
$x_{n_k}\in\D_\alpha^{n_k}$ and $\delta(x_{n_k},\D_\alpha)>\epsilon$ for all
$k\in\IN$. Since the sequence $(x_{n_k})_{k\in\IN}$ is bounded there is a
convergent subsequence. Again, we assume w.l.o.g. that the sequence itself is
convergent to a point $x_0$. From the continuity of the mapping $x\mapsto
\delta(x,\D_\alpha)$ it follows that $\delta(x_0,\D_\alpha)\ge
\epsilon$. Since the interior of $\D_\alpha$ is non-empty, we can choose $d$
points $z_1,\dots,z_d$ in $\interior\D_\alpha$ in such a way that $x_0$ and
$z_1,\dots,z_d$ are in general position. It is easy to show that
$\lim_{k\to\infty}x_{n_k}=x_0$ implies that
\[
\Hlim_{k\to\infty}S[x_{n_k},z_1,\dots,z_d]=S[x_0,z_1,\dots,z_d]
\]
where $S[x_0,z_1,\dots,z_d]$ denotes the simplex generated by the points $x_0,
z_1,\dots,\linebreak z_d$. Since $\D(z_i)>\alpha$ for $i=1,\dots,d$, \PtwD implies that
there exists $N$ such that $\D^n(z_i)\ge\alpha$ for $i=1,\dots,d$ and $n\ge N$.
Because of the convexity of the trimmed regions we also have
$S[x_{n_k},z_1,\dots,z_d]\subset\D_\alpha^{n_k}$ for all $k$ with $n_k\ge N$.
Now, let $z_0\in\interior S[x_0,z_1,\dots,z_d]\setminus\D_\alpha$. Because of
the Hausdorff convergence of the simplices, it follows from Corollary~\ref{cor.HD}
that there is a $K\in\IN$ such that $z_0\in\D_\alpha^{n_k}$ for $k\ge K$.
But then, $z_0\in\limsup_{n\to\infty}\D_\alpha^n$ in contradiction to
$z_0\notin\D_\alpha$. Thus, \eqref{eq.konv5} has to be valid and the Hausdorff
convergence of the trimmed regions for $\alpha\in(0,\alpha_{\max})$ is shown.

To show the second part of the theorem we assume w.l.o.g. that
$\D_{\alpha_{\max}}^n\ne\emptyset$ for all $n\in\IN$. According to
Theorem~\ref{th.contalpha} the mapping $\alpha\mapsto\D_\alpha$ is left
continuous on $(0,\alpha_{\max}]$. Thus, there is an $\epsilon>0$ and
$\beta<\alpha_{\max}$, such that $\delta_H(\D_\beta,\D_{\alpha_{\max}})<
\frac\epsilon6$. From what has already been proven there exists $N\in\IN$,
such that $\delta_H(\D_\beta^n,\D_\beta)<\frac\epsilon6$ for all $n\ge N$.
Because of $\D_{\alpha_{\max}}^n\subset\D_\beta^n$ we have
$\delta_H(\D_{\alpha_{\max}}^n,\D_\beta^n)\le\diam(\D_\beta^n)$.
It is easy to show that $\delta_H(A,B)<\epsilon$ implies that
$\diam(A)<\diam(B)+2\epsilon$. For $n\ge N$ we thus get
\[
\delta_H(\D_{\alpha_{\max}}^n,\D_\beta^n)\le\diam(\D_\beta^n)
<\diam(\D_\beta)+\frac13\epsilon
<\diam(\D_{\alpha_{\max}})+\frac23\epsilon
=\frac23\epsilon\,,
\]
since the diameter of a singleton is equal to zero. Therefore,
\[
\delta_H(\D_{\alpha_{\max}},\D_{\alpha_{\max}}^n)
\le
\underbrace{\delta_H(\D_{\alpha_{\max}},\D_\beta)}_{<\frac\epsilon6} +
\underbrace{\delta_H(D_\beta,D_\beta^n)}_{<\frac\epsilon6} +
\underbrace{\delta_H(D_\beta^n,\D_{\alpha_{\max}}^n)}_{<\frac{4\epsilon}6}
< \epsilon
\]
for all $n\ge N$ and $\Hlim_{n\to\infty}\D_{\alpha_{\max}}^n=\D_{\alpha_{\max}}$
is shown.\qed

\textbf{Proof (of Theorem~\ref{th.PktZPktT}):} Because of Theorem~\ref{th.konvpunkt}
we just have to show that for every $\alpha>0$ holds:
\[
\{z\in\IR^d\,|\,\D(z)>\alpha\}\subset\liminf_{n\to\infty}\D_\alpha^n
\subset\limsup_{n\to\infty}\D_\alpha^n\subset\D_\alpha\,.
\]
We start to show $\{z\in\IR^d\,|\,\D(z)>\alpha\}\subset\liminf_{n\to\infty}\D_\alpha^n$.
If $\alpha\ge\alpha_{\max}$, then $\{z\in\IR^d\,|\,\D(z)>\alpha\}$ is empty
and the assertion is trivially satisfied. If $\alpha<\alpha_{\max}$ and
$\D(z)>\alpha$ then because of the continuity of $\D(\,\cdot\,|\,P)$ the point
$z$ lies in the interior of $\D_\alpha$. By assumption $\Hlim_{n\to\infty}
\D_\alpha^n=\D_\alpha$. From Corollary~\ref{cor.HD} then follows that there
is an $N$ such that $z\in\D_\alpha^n$ for all $n\ge N$. But then,
$z\in\liminf\D_\alpha^n$ as well and
$\{z\in\IR^d\,|\,\D(z)>\alpha\}\subset\liminf_{n\to\infty}\D_\alpha^n$.

We next show $\limsup_{n\to\infty}\D_\alpha^n\subset\D_\alpha$.
If $\alpha>\alpha_{max}$, then $\D_\alpha^n=\emptyset$ for $n\ge N_\alpha$.
Thus, $\limsup\D_\alpha^n=\emptyset$ and the assertion is satisfied.
Now, let $\alpha\le\alpha_{\max}$ and $\D(z)<\alpha$. Then there is $\beta$
with $\D(z)<\beta<\alpha$. Thus, $z$ lies in the complement of $\D_\beta$.
Again, it follows from Corollary~\ref{cor.HD} that there is an $N$, such that
$z\notin\D_\beta^n$ for all $n\ge N$. Thus, $z\notin\limsup\D_\beta^n$.
Since $\limsup\D_\alpha^n\subset\limsup\D_\beta^n$ it follows that
$z\notin\limsup\D_\alpha^n$ and the assertion is proved. \qed

\textbf{Proof (of Theorem~\ref{th.konvglm}):} $(i)\implies (ii)$ is trivial.

We show $(ii)\implies (iii)$. Since every non-trivial depth assumes at least
two values, there is $x_0\in \IR^d$ with $\alpha_0:=\D(x_0)>0$. We show the
assertion w.l.o.g. for $0<\alpha<\alpha_0$ and $\epsilon$ such that
$\alpha-\epsilon>0$ and $\alpha+\epsilon<\alpha_0$. In that case
$\D_{\alpha-\epsilon}$ is bounded and there is a compact set $M$ such that
$\D_{\alpha-\epsilon}$ is contained in the interior of $M$.
Thus, there is $N_\epsilon$, such that
\[
\sup_{x\in M}|\D^n(x)-\D(x)|\le\epsilon\quad\text{for all $n\ge N_\epsilon$.}
\]
Now let $n\ge N_\epsilon$ and $x\in\D_{\beta+\epsilon}$ with $\beta\ge\alpha$.
Then, $x\in M$ and it holds
\[
|\D^n(x)-\D(x)|\le\epsilon\,.
\]
In particular,
\[
\D^n(x)\ge\D(x)-\epsilon\ge(\beta+\epsilon)-\epsilon=\beta\,,
\]
i.e., $x\in\D_\beta^n$. Therefore, it is shown that
$\D_{\beta+\epsilon}\subset D_\beta^n$ for all $n\ge N_\epsilon$ and
$\beta\ge\alpha$.

In the following let $n\ge N_\epsilon$. To show that $\D_\beta^n\subset
\D_{\beta-\epsilon}$ for all $\beta\ge\alpha$, first note that
\[
\D^n(x)\le\D(x)+\epsilon\quad\text{for all $x\in M$.}
\]
For $x\in M\setminus\D_{\alpha-\epsilon}$ then holds
\[
\D^n(x)\le\D(x)+\epsilon<(\alpha-\epsilon)+\epsilon=\alpha
\]
and therefore $M\setminus\D_{\alpha-\epsilon}\subset(\D_\alpha^n)^c$. The trimmed
regions are star-shaped and therefore connected. Thus, $\D_\alpha^n$ is
either a subset of $\D_{\alpha-\epsilon}$ or a subset of $M^c$.
Assume that $\D_\alpha^n\subset M^c$. From the choice of $\alpha$ it is clear
that $x_0\in\D_{\alpha+\epsilon}$ and because of $\D_{\alpha+\epsilon}\subset
\D_\alpha^n$ we have
\[
x_0\in\D_{\alpha+\epsilon}\subset\D_\alpha^n\subset M^c\subset(\D_{\alpha-\epsilon})^c\,,
\]
i.e., $\D(x_0)<\alpha-\epsilon<\alpha_0$, contradiction. Thus,
$D_\alpha^n\subset\D_{\alpha-\epsilon}$ for all $n\ge N_\epsilon$.

Now let $\beta\ge\alpha$. Then, $\D_\beta^n\subset\D_\alpha^n\subset
\D_{\alpha-\epsilon}\subset M$. If $x\in\D_\beta^n$, then
\[
\D(x)\ge\D^n(x)-\epsilon\ge\beta-\epsilon\,,
\]
i.e., $x\in\D_{\beta-\epsilon}$. Thus, it is shown that $\D_\beta^n\subset
\D_{\beta-\epsilon}$ for all $n\ge N_\epsilon$ and $\beta\ge\alpha$.
All in all we get
\[
\D_{\beta+\epsilon}\subset\D_\beta^n\subset\D_{\beta-\epsilon}
\quad\text{for all $n\ge N_\epsilon$ and $\beta\ge\alpha$,}
\]
as stated.

$(iii)\implies (iv)$ is again trivial.

We now show $(iv)\implies(i)$. We have to show that for every $\epsilon>0$
there exists $N_\epsilon$, such that
\[
\sup_{x\in\IR^d}|\,\D^n(x)-\D(x)\,|\le\epsilon\quad\text{for all $n\ge N_\epsilon$.}
\]
Let $\epsilon>0$ be given. We choose $A=[\epsilon/4,\alpha_{\max}-\epsilon/4]$.
According to $(iv)$ there exists $N\in\IN$, such that
\begin{equation}\label{eq.konv3}
\D_{\beta+\frac\epsilon4}\subset\D_\beta^n\subset\D_{\beta-\frac\epsilon4}
\quad\text{for all $n\ge N$ and $\beta\in A$,}
\end{equation}
and $\D_{\alpha_{\max}+\frac\epsilon4}^n=\emptyset$ for all $n\ge N$. In the 
following let $n\ge N$.

\textbf{Case 1:} If $x\in\D_{\frac\epsilon2}$, then $\D(x)=:\gamma\ge
\frac\epsilon2$. Thus, $x\in\D_\gamma$ and because of the assumption
also $x\in\D_{\gamma-\frac\epsilon4}^n$. Therefore,
\[
\D^n(x)\ge\gamma-\frac\epsilon4=\D(x)-\frac\epsilon4\,.
\]
To bound $\D^n(x)$ from above we distinguish two cases:

\textbf{Case 1a:} If $\gamma<\alpha_{\max}-\frac34\epsilon$, it follows from
$x\notin\D_{\gamma+\frac\epsilon4}$, that $x\notin\D_{\gamma+\frac\epsilon2}^n$.
From this follows
\[
\D^n(x)<\gamma+\frac\epsilon2=\D(x)+\frac\epsilon2\,.
\]

\textbf{Case 1b:} If $\gamma\ge\alpha_{\max}-\frac34\epsilon$, then
$\gamma+\epsilon\ge\alpha_{\max}+\frac\epsilon4$ and thus
$\D_{\gamma+\epsilon}^n\subset\D_{\alpha_{\max}+\frac\epsilon4}^n=\emptyset$.
It follows that
\[
\D^n(x)<\gamma+\epsilon=\D(x)+\epsilon\,.
\]
From Case~1 together with Cases~1a and 1b it follows
\[
\D(x)-\frac\epsilon4\le\D^n(x)<\D(x)+\epsilon\quad\text{for all $x\in\D_{\frac\epsilon2}$}
\]
and therefore also
\[
\sup_{x\in\D_{\frac\epsilon2}}|\,\D^n(x)-\D(x)\,|\le\epsilon\,.
\]

\textbf{Case 2:} If $x\notin\D_{\frac\epsilon2}$, then $\D(x)<\frac\epsilon2$.
Further, with $\beta=\epsilon/4$ it follows from equation~\eqref{eq.konv3} that 
$x\notin\D_{\frac{3\epsilon}4}^n$, i.e., $\D^n(x)<\frac{3\epsilon}4$. From this 
we conclude
\[
|\,\D^n(x)-\D(x)\,|<\frac{3\epsilon}{4}\quad\text{for all $x\notin\D_{\frac\epsilon2}$.}
\]
and
\[
\sup_{x\notin\D_{\frac\epsilon2}}|\,\D^n(x)-\D(x)\,|\le\frac{3\epsilon}{4}\,.
\]
From the two cases we finally get
\[
\sup_{x\in\IR^d}|\,\D^n(x)-\D(x)\,|\le\epsilon
\quad\text{for all $n\ge N$,}
\]
as was to be shown.\qed

\textbf{Proof (of Theorem~\ref{th.KomTKomZ}):}
Let $[\alpha_1,\alpha_2]\subset(0,\alpha_{\max})$ and $\epsilon>0$.
Because of the strict monotonicity the mapping $\alpha\mapsto\D_\alpha$ is
continuous. Since every continuous function on a compact set is uniformly
continuous, this mapping is uniformly continuous on $[\alpha_1/2,
\alpha_{\max}]$. Thus, there exists $\gamma>0$, such that
\[
\delta_H(\D_{\beta_1},\D_{\beta_2})<\frac\epsilon2
\quad\text{for all $\beta_1,\beta_2\in[\alpha_1/2,\alpha_{\max}]$
with $|\beta_1-\beta_2|\le2\gamma$.}
\]
Assume w.l.o.g. that $\gamma$ is so small, that $\gamma<\alpha_1/2$ and 
$\alpha_2\le\alpha_{\max}-\gamma$. 
Then it follows that
\[
\delta_H(\D_{\beta-\gamma},\D_{\beta+\gamma})<\frac\epsilon2
\quad\text{for all $\beta\in[\alpha_1,\alpha_{\max}-\gamma]$.}
\]
If \ComD is satisfied then it follows from Theorem~\ref{th.konvglm}, 
Part~$(iii)$, that there exists $N\in\IN$, such that
\[
\D_{\beta+\gamma}\subset\D_\beta^n\subset\D_{\beta-\gamma}
\quad\text{for all $n\ge N$ and $\beta\ge\alpha_1$.}
\]
Trivially, $\D_{\beta+\gamma}\subset\D_\beta\subset\D_{\beta-\gamma}$
holds as well. Thus, for every $n\ge N$ and $\beta\in[\alpha_1,\alpha_{\max}-\gamma]$:
\[
\delta_H(\D_\beta^n,\D_\beta)
\le\delta_H(\D_{\beta+\gamma},\D_{\beta-\gamma})
<\frac\epsilon2\,.
\]
Therefore,
\[
\sup_{\alpha_1\le\beta\le\alpha_2}\delta_H(\D_\beta^n,\D_\beta)\le\frac\epsilon2
\quad\text{for all $n\ge N$,}
\]
which implies \ComR.

For proving the second part of the Theorem, note that $\D_{\beta+\gamma}$ will 
be empty, when $\beta>\alpha_{\max}-\gamma$. Thus, for 
$\beta\in(\alpha_{\max-\gamma},\alpha_{\max}]$ we only have
\[
\D_\beta^n\subset\D_{\beta-\gamma}\quad\text{for all $n\ge N$.}
\]
Trivially, $\D_\beta\subset\D_{\beta-\gamma}$ holds as well.
Now, if $\D_{\alpha_{\max}}$ is a singleton, i.e., $\D_{\alpha_{\max}}=\{x_0\}$, then
\begin{align*}
\delta_H(\D_\beta^n,\D_\beta)
&\le
\delta_H(\D_\beta^n,\D_{\alpha_{\max}})+\delta_H(\D_{\alpha_{\max}},\D_\beta)
<
\delta_H(\D_\beta^n,\{x_0\})+\frac\epsilon2\\
&=
\max_{x\in\D_\beta^n}\|x-x_0\|+\frac\epsilon2
\le
\max_{x\in\D_{\beta-\gamma}}\|x-x_0\|+\frac\epsilon2\\
&=
\delta_H(\D_{\beta-\gamma},\{x_0\})+\frac\epsilon2
<\frac\epsilon2+\frac\epsilon2
=\epsilon\,,
\end{align*}
and the proof is finished.\qed

\textbf{Proof (of Theorem~\ref{th.KomZGlmT}):} We show that \ComR together
with \RC implies Condition $(iv)$
in Theorem~\ref{th.konvglm}. Thus, we have to show that for each $\epsilon>0$
and every compact interval $A\subset(0,\alpha_{\max})$ there exists
$N_{A,\epsilon}\in\IN$, such that
\[
\D_{\alpha+\epsilon}\subset\D_\alpha^n\subset\D_{\alpha-\epsilon}
\quad\text{for all $n\ge N_{A,\epsilon}$ and $\alpha\in A$.}
\]
Let $A=[\alpha_1,\alpha_2]$ and $\epsilon$ be given. We assume w.l.o.g.
that $\epsilon<\alpha_1$ and $\epsilon<\alpha_{\max}-\alpha_2$.
Since a continuous function, defined on a compact set, is uniformly continuous
and since $\D_{\alpha_1-\epsilon}$ is compact, the mapping
$z\mapsto\D(z\,|\,P)$ is uniformly continuous on $\D_{\alpha_1-\epsilon}$.
Thus, there is $\gamma>0$, such that
\[
|\D(x)-\D(y)|<\epsilon\quad\text{for all $x,y\in\D_{\alpha_1-\epsilon}$ with
$\|x-y\|<\delta$.}
\]
Further,
\[
\min_{x\in\partial\D_{\alpha-\epsilon},y\in\partial\D_\alpha}\|x-y\|\ge\gamma
\quad\text{for all $\alpha\in[\alpha_1,\alpha_{\max}]$.}
\]
This holds because if the above equation was not satisfied then there was
$x\in\partial\D_{\alpha-\epsilon}$ and $y\in\partial\D_\alpha$ such that
$\|x-y\|<\gamma$. Since $\D(\,\cdot\,|\,P)$ is continuous this would imply
$\D(x)=\alpha-\epsilon$ and $\D(y)=\alpha$. Therefore we would get
\[
|\D(x)-\D(y)|=\D(y)-\D(x)=\alpha-(\alpha-\epsilon)=\epsilon
\]
in contradiction to the uniform continuity.

Now, choose $N_{A,\epsilon}$ so large that $\delta_H(\D_\alpha^n,\D_\alpha)
<\gamma$ for all $n\ge N_{A,\epsilon}$ and $\alpha\in A$. Since the trimmed
regions are convex it follows from Proposition~\ref{th.HD-stetigkeit} in the
Appendix that
\[
\D_{\alpha+\epsilon}\subset\D_\alpha^n\subset\D_{\alpha-\epsilon}
\quad\text{for all $n\ge N_{A,\epsilon}$ and $\alpha\in A$,}
\]
as was to be shown.\qed

In the proofs of Theorems~\ref{th.PktZKomZ} and \ref{th.PktTGlmT} we make use of
another notion of convergence, the so-called \emph{continuous convergence}.
\begin{definition}
Let $(X,\rho_X)$ and $(Y,\rho_Y)$ be metric spaces. A sequence $(f_n)_{n\in\IN}$
of mappings from $X$ to $Y$ is said to \emph{converge continuously} to $f$ if
for each $x\in X$ and for each sequence $(x_n)_{n\in\IN}$ such that
$\lim_{n\to\infty}x_n=x$ we have $\lim_{n\to\infty}f_n(x_n)=f(x)$.
\end{definition}
The following well-known result that connects continuous convergence and compact 
convergence will be useful in the proofs of Theorems~\ref{th.PktZKomZ} and 
\ref{th.PktTGlmT}.
\begin{proposition}\label{pr.contconv}
Let $(X,\rho_X)$ and $(Y,\rho_Y)$ be metric spaces and $(f_n)_{n\in\IN}$
be a sequence of functions from $X$ to $Y$.
\begin{enumerate}
\item If $f_n$ converges continuously to $f$, then $f$ is continuous.
\item If the sequence $(f_n)$ converges continuously to $f$, then
$(f_n)$ is compact convergent to $f$.
\item Let $X$ be locally compact. If the sequence $(f_n)$ is compact convergent 
to $f$ and $f$ is continuous, then $(f_n)$ converges continuously to $f$.  
\end{enumerate}
\end{proposition}

\textbf{Proof (of Theorem~\ref{th.PktZKomZ}):} We show that under the given
assumptions the sequence
of mappings $\alpha\mapsto\D^n_\alpha$,$n\in\IN$, is continuous convergent to
the mapping $\alpha\mapsto\D_\alpha$, i.e., for every sequence $(\alpha_n)$ that
is convergent to $\alpha_0$ it holds that $\Hlim_{n\to\infty}\D^n_{\alpha_n}=
\D_\alpha$. From Proposition~\ref{pr.contconv} above it then follows that the
trimmed regions are compact convergent.

First assume that $0<\alpha_0<\alpha_{\max}$. Because of the strict monotonicity
it follows from Theorem~\ref{th.contalpha} that the mapping $\alpha\mapsto
\D_\alpha$ is continuous. Thus, there is $\gamma>0$ such that
$\delta_H(\D_{\alpha_0-\gamma}, \D_{\alpha_0+\gamma})<\epsilon/5$. Further, from
the pointwise convergence of the trimmed regions it follows that there is an
$N_1\in\IN$ such that
$\delta_H(\D_{\alpha_0-\gamma},\D_{\alpha_0-\gamma}^n)<\epsilon/5$ and
$\delta_H(\D_{\alpha_0+\gamma},\D_{\alpha_0+\gamma}^n)<\epsilon/5$ for all
$n\ge N_1$. We conclude that for $n\ge N_1$
\begin{multline*}
\delta_H(\D_{\alpha_0-\gamma}^n,\D_{\alpha_0+\gamma}^n)
\le
\delta_H(\D_{\alpha_0-\gamma}^n,\D_{\alpha_0-\gamma})\\
+
\delta_H(\D_{\alpha_0-\gamma},\D_{\alpha_0+\gamma})
+
\delta_H(\D_{\alpha_0+\gamma},\D_{\alpha_0+\gamma}^n)
<\frac{3\epsilon}5\,.
\end{multline*}
Since $(\alpha_n)$ converges to $\alpha_0$ there is an $N_2\in\IN$ such that
$|\alpha_n-\alpha_0|<\gamma$. For $n\ge N_2$ the trimmed region
$\D_{\alpha_n}^n$ lies between $\D_{\alpha_0-\gamma}^n$ and
$\D_{\alpha_0+\gamma}^n$. Therefore
\[
\delta_H(\D_{\alpha_n}^n,\D_{\alpha_0-\gamma}^n)\le
\delta_H(\D_{\alpha_0+\gamma}^n,\D_{\alpha_0-\gamma}^n)\le
\frac{3\epsilon}5
\]
for $n\ge N_2$. From this follows that for $n\ge\max\{N_1,N_2\}$
\[
\delta_H(\D_{\alpha_n}^n,\D_{\alpha_0})
\le
\delta_H(\D_{\alpha_n}^n,\D_{\alpha_0-\gamma}^n)
+
\delta_H(\D_{\alpha_0-\gamma}^n,\D_{\alpha_0-\gamma})+
\delta_H(\D_{\alpha_0-\gamma},\D_{\alpha_0})
<\epsilon\,.
\]
Now assume that $\alpha_0=\alpha_{\max}$ and $\D_{\alpha_0}=\{x_0\}$ is
a singleton. As above there is $\gamma>0$ and $N_1\in\IN$ such that
$\delta_H(\D_{\alpha_0},\D_{\alpha_0-\gamma})<\epsilon/2$ and
$\delta_H(\D_{\alpha_0-\gamma},\D_{\alpha_0-\gamma}^n)<\epsilon/2$
for $n\ge N_1$. Choose $N_2\in\IN$ such that $\alpha_n>\alpha_0-\gamma$.
Then, for $n\ge N_2$ the trimmed region $\D_{\alpha_n}^n$ is contained in
$\D_{\alpha_0-\gamma}^n$. For $n\ge\max\{N_1,N_2\}$ we thus get
\[
\delta_H(\D_{\alpha_n}^n,\D_{\alpha_0})
\le
\delta_H(\D_{\alpha_0-\gamma}^n,\D_{\alpha_0})
\le
\delta_H(\D_{\alpha_0-\gamma}^n,\D_{\alpha_0-\gamma})+
\delta_H(\D_{\alpha_0-\gamma},\D_{\alpha_0})
<\epsilon
\]
and the second part of the theorem is proved.\qed

\textbf{Proof (of Theorem~\ref{th.PktTGlmT}):} We show that under the assumptions
the sequence
$(\D^n)_{n\in\IN}$ is continuous convergent to $\D$, i.e., for every sequence
$(z_n)$ that is convergent to $z_0$ it holds that $\lim_{n\to\infty}\D^n(z_n)
=\D(z_0)$. From Proposition~\ref{pr.contconv} above it then follows that
$(\D^n)$ is compact and thus uniform convergent to $\D$.

We start with showing that for every $\epsilon>0$ there is an $N\in\IN$ such
that  $D^n(z_n)\ge\D(z_0)-\epsilon$ for all $n\ge N$. Let $\D(z_0)=\alpha$.
Since $\D$ is continuous, $z_0\in\interior\D_{\alpha-\epsilon/4}$. Then, there
is $\gamma>0$ such that $B(z_0,\gamma)\subset\interior\D_{\alpha-\epsilon/2}$.
Choose $d+1$ points $x_1,\dots,x_{d+1}\in B(z_0,\gamma)$ such that $z_0\in
\interior S[x_1,\dots,x_{d+1}]$, where $S[x_1,\dots,x_{d+1}]$ denotes the
simplex generated by $x_1,\dots,x_{d+1}$.

From $\D(x_i)\ge\alpha-\frac\epsilon2$, $i=1,\dots,d+1$, and the pointwise
convergence of the depths it follows that there exists $N_1\in\IN$ such that
$\D^n(x_i)\ge\alpha-\epsilon$ for all $n\ge N_1$ and $i=1,\dots,d+1$.
For $n\ge N_1$ we thus have $x_i\in\D^n_{\alpha-\frac\epsilon2}$, $i=1,\dots,d+1$.
Since the trimmed regions of $\D$ are convex it follows that
$S[x_1,\dots,x_{d+1}]\subset\D^n_{\alpha-\epsilon}$ for $n\ge N_1$.
Now, because of $z_0\in\interior S[x_1,\dots,x_{d+1}]$ there is $N_2\in\IN$,
such that $z_n\in S[x_1,\dots,x_{d+1}]$ for all $n\ge N_2$. For $n\ge\max\{N_1,N_2\}$
we therefore have $z_n\in \D^n_{\alpha-\epsilon}$, i.e.,
$\D^n(z_n)\ge\D(z_0)-\epsilon$.

It remains to show that there exists $N\in\IN$ such that
$D^n(z_n)\le\D(z_0)+\epsilon$ for all $n\ge N$. Assume the contrary. Then there
is a subsequence $(n_k)_{k\in\IN}$, such that
\[
\D^{n_k}(z_{n_k})>\D(z_0)+\epsilon=\alpha+\epsilon\quad\text{for all $k\in\IN$.}
\]
For $\alpha=\alpha_{\max}$ this is a contradiction to \RC.
If $\alpha<\alpha_{\max}$ we assume w.l.o.g. that $\alpha+\epsilon<
\alpha_{\max}$. Since the trimmed regions are closed, their complements are
open. Thus, there is $\gamma>0$ such that $B(z_0,\gamma)\intersection
\D_{\alpha+\epsilon}=\emptyset$. For $\beta$ with $\alpha+\epsilon<\beta<
\alpha_{\max}$ the trimmed region $\D_\beta$ is larger than $\D_{\alpha_{\max}}$
and the interior of $\D_\beta$ is non-empty. Choose $d$ points $x_1,\dots,x_d\in
\D_\beta$ such that $z_0,x_1,\dots,x_d$ are in general position.
Because of pointwise convergence there exists $K_1\in\IN$ such that
$\D^{n_k}(x_i)>\alpha+\epsilon$ for $k\ge K_1$ and $i=1,\dots,d$. From the
convexity of the trimmed regions it follows that $S[z_{n_k},x_1,\dots,x_d]
\subset\D_{\alpha+\epsilon}^{n_k}$ for $k\ge K_1$. It is easy to show that
\[
\Hlim_{k\to\infty}S[z_{n_k},x_1,\dots,x_d]=S[z_0,x_1,\dots,x_d]\,.
\]
Since $\interior S[z_0,x_1,\dots,x_d]\intersection B(z_0,\gamma)$ is non-empty
there exists $z^*\in\linebreak \interior S[z_0,x_1,\dots,x_d]\intersection
B(z_0,\gamma)$. Because of the Hausdorff convergence of the simplices
it follows from Corollary~\ref{cor.HD} that there is $K_2\in\IN$ such that
$z^*\in S[z_{n_k},x_1,\dots,x_d]$ for all $k\ge K_2$. Thus,
$\D^{n_k}(z^*)\ge\alpha+\epsilon$ for $k\ge K=\max\{K_1,K_2\}$. From
pointwise convergence of the depth it follows $\D(z^*)\ge\alpha+\epsilon$ in 
contradiction to $B(z_0,\gamma)\intersection\D_{\alpha+\frac\epsilon2}=\emptyset$. 
Thus, the proof is finished.\qed

\section{Hausdorff-convergence}\label{sec.Hausdorff}
In this section we state the definition of Hausdorff convergence as well as some
important facts on this notion of convergence. Detailed studies of the
notion of Hausdorff convergence can be found, e.g., in \cite{KleinT84}
and \cite{Beer93}.

The Euclidean distance between two points $x,y\in\IR^d$ is given by
$\delta(x,y)=\|x-y\|$. The distance between a point $x\in\IR^d$ and a set
$A\subset\IR^d$ can then be defined by $\delta(x,A)=\inf_{y\in A}\delta(x,y)$.
If $A$ is closed one can write $\min$ instead of $\inf$. The set of all
non-empty compact subsets of $\IR^d$ is denoted by ${\cal K}_0^d$.
\begin{definition}[Hausdorff distance]\label{def.HD-Distanz}%
For $A,B\in{\cal K}_0^d$ the \emph{Hausdorff distance} $\delta_H(A,B)$ is defined
by
\[
\delta_H(A,B)=\max\{\max_{x\in A}\delta(x,B),\max_{x\in B}\delta(x,A)\}\,.
\]
\end{definition}
For $A\subset\IR^d$ and $\epsilon>0$ the \emph{$\epsilon$-neighborhood}
$U_\epsilon(A)$ is given by $U_\epsilon(A)=\{x\in\IR^d\,|\,\delta(x,A)<\epsilon\}
=\Union_{x\in A}B(x,\epsilon)$. With this notation an equivalent definition of
the Hausdorff distance is $\delta_H(A,B) = \inf\{\,\epsilon>0\,|\,A\subset U_\epsilon(B),
B\subset U_\epsilon(A)\}$.

If one of the sets is a singleton, then $\delta_H(\{x_0\},A)=
\max_{x\in A}\delta(x,x_0)$.

If $A_1\subset B\subset A_2$ then $\delta_H(A_1,B)\le\delta_H(A_1,A_2)$
as well as $\delta_H(B,A_2)\le\delta_H(A_1,A_2)$.
If $A_1\subset B_i\subset A_2$, $i=1,2$, then $\delta_H(B_1,B_2)\le
\delta_H(A_1,A_2)$.

The Hausdorff distance is a metric on ${\cal K}_0^d$. Thus, the pair
$({\cal K}_0^d,\delta_H)$ is a metric space. Therefore it is possible to define
convergence of compact sets in the Hausdorff metric or short \emph{Hausdorff
convergence}.
\begin{definition}[Hausdorff convergence]\label{def.HD-Konvergenz}
Let $(K_n)_{n\in\IN}$ be a sequence of non-empty compact subsets of $\IR^d$.
The sequence $(K_n)_{n\in\IN}$ is said to be \emph{Hausdorff convergent}
to a set $K\in{\cal K}_0^d$, if $\lim_{n\to\infty}\delta_H(K_n,K)=0$. In this
case we write $\Hlim_{n\to\infty}K_n=K$.
\end{definition}

\begin{proposition}\label{pr.HD-limes-aufab}
If a sequence $(K_n)_{n\in\IN}$ is decreasing, i.e., $K_1\supset K_2\supset\dots$,
then the Hausdorff limit exists and is given by
$\Hlim_{n\to\infty}K_n=\Intersection_{n=1}^\infty K_n$.

If a sequence $(K_n)_{n\in\IN}$ is increasing, i.e.,
$K_1\subset K_2\subset\dots$, and the union of the sets is bounded, then the
Hausdorff limit exists and is given by $\Hlim_{n\to\infty}K_n=
\cl(\Union_{n=1}^\infty K_n)$.
\end{proposition}

If the sets $K_n$, $n\in\IN$, are connected then the following criteria is
useful.
\begin{proposition}\label{th.AWHD}
Let $(K_n)_{n\in\IN}$ be sequence of connected sets in ${\cal K}_0^d$
and let $K\in {\cal K}_0^d$. If for each $M>0$ it holds that
\[
\lim_{n\to\infty}\max\{\delta(x,K_n)\,|\,x\in K\intersection B(0,M)\}=0
\]
and
\[
\lim_{n\to\infty}\max\{\delta(x,K)\,|\,x\in K_n\intersection B(0,M)\}=0\,,
\]
then $\Hlim_{n\to\infty}K_n=K$.
\end{proposition}

The following proposition and its corollary show that for convex sets
Hausdorff convergence behaves nicely.
\begin{proposition}\label{th.HD-stetigkeit}
If $A_1,A_2$ are convex sets in ${\cal K}_0^d$ such that $A_1\subset A_2$ and
\[
\min_{x\in\partial A_1,y\in\partial A_2}\delta(x,y)=\gamma>0\,,
\]
then for every convex set $B\in{\cal K}_0^d$ holds
\begin{align*}
\delta_H(B,A_1)\le \gamma&\implies B\subset A_2\,,\\
\delta_H(B,A_2)\le \gamma&\implies A_1\subset B\,.
\end{align*}
\end{proposition}

\begin{corollary}\label{cor.HD}
Let $(K_n)_{n\in\IN}$ be a sequence  of convex sets in ${\cal K}_0^d$ with
$\Hlim_{n\to\infty}K_n=K$, where $K$ is a convex set in ${\cal K}_0^d$.
Then the following assertions hold:
\begin{enumerate}
\item For every $x\in\interior K$ there is an $N\in\IN$, such that
$x\in K_n$ for all $n\ge N$.
\item For every $x\in K^c$ there is an $N\in\IN$, such that
$x\notin K_n$ for all $n\ge \IN$.
\item $\interior K\subset\liminf_{n\to\infty} K_n\subset\limsup_{n\to\infty} K_n\subset K$.
\end{enumerate}
\end{corollary}


\newpage
\section*{Supplementary material}
\textbf{Supplement to ``Convergence of depths and depth-trimmed regions''}
The supplement contains some examples that show that
without the assumption `strictly monotone for $P$' Theorems~\ref{th.PktTPktZ}, 
\ref{th.KomTKomZ}, and \ref{th.PktZKomZ} are in general false, and 
without the assumption `continuous for $P$' Theorems~\ref{th.PktZPktT},  
\ref{th.KomZGlmT}, and  \ref{th.PktTGlmT} are in general false.
\end{appendix}
\manuallabel{ex.1}{2.1}
\manuallabel{ex.2}{2.2}
\manuallabel{ex.3}{2.3}
\manuallabel{ex.4}{2.4}


\bibliographystyle{imsart-nameyear}
\bibliography{CDaDTR}

\begin{thebibliography}{46}

\bibitem[\protect\citeauthoryear{Barnett}{1976}]{Barnett76}
\begin{barticle}[author]
\bauthor{\bsnm{Barnett},~\bfnm{V.}\binits{V.}}
(\byear{1976}).
\btitle{The ordering of multivariate data}.
\bjournal{Journal of the Royal Statistical Society, Series A}
\bvolume{139}
\bpages{318--352}.
\bnote{With Discussion}.
\end{barticle}
\endbibitem

\bibitem[\protect\citeauthoryear{Beer}{1993}]{Beer93}
\begin{bbook}[author]
\bauthor{\bsnm{Beer},~\bfnm{Gerald}\binits{G.}}
(\byear{1993}).
\btitle{Topologies on Closed and Closed Convex Sets}.
\bpublisher{Kluwer}, \baddress{Dordrecht}.
\end{bbook}
\endbibitem

\bibitem[\protect\citeauthoryear{Cascos and
  L{\'{o}}pez-D{\'{\i}}az}{2008}]{CascosLD08}
\begin{barticle}[author]
\bauthor{\bsnm{Cascos},~\bfnm{Ignacio}\binits{I.}} \AND
  \bauthor{\bsnm{L{\'{o}}pez-D{\'{\i}}az},~\bfnm{Miguel}\binits{M.}}
(\byear{2008}).
\btitle{Consistency of the $\alpha$-trimming of a probability. Applications to
  central regions}.
\bjournal{Bernoulli}
\bvolume{14}
\bpages{580--592}.
\end{barticle}
\endbibitem

\bibitem[\protect\citeauthoryear{Cascos and
  L{\'{o}}pez-D{\'{\i}}az}{2016}]{CascosLD16}
\begin{barticle}[author]
\bauthor{\bsnm{Cascos},~\bfnm{Ignacio}\binits{I.}} \AND
  \bauthor{\bsnm{L{\'{o}}pez-D{\'{\i}}az},~\bfnm{Miguel}\binits{M.}}
(\byear{2016}).
\btitle{On the uniform consistency of the zonoid depth}.
\bjournal{Journal of Multivariate Analysis}
\bvolume{143}
\bpages{394--397}.
\end{barticle}
\endbibitem

\bibitem[\protect\citeauthoryear{Cascos and Molchanov}{2007}]{CascosM07}
\begin{barticle}[author]
\bauthor{\bsnm{Cascos},~\bfnm{Ignacio}\binits{I.}} \AND
  \bauthor{\bsnm{Molchanov},~\bfnm{Ilya}\binits{I.}}
(\byear{2007}).
\btitle{Multivariate risks and depth-trimmed regions}.
\bjournal{Finance and Stochastics}
\bvolume{11}
\bpages{373--397}.
\end{barticle}
\endbibitem

\bibitem[\protect\citeauthoryear{Cramer}{2003}]{Cramer03}
\begin{bbook}[author]
\bauthor{\bsnm{Cramer},~\bfnm{Katharina}\binits{K.}}
(\byear{2003}).
\btitle{Multivariate Ausrei{\ss}er und Datentiefe}.
\bpublisher{Shaker}, \baddress{Aachen}.
\end{bbook}
\endbibitem

\bibitem[\protect\citeauthoryear{Donoho}{1982}]{Donoho82}
\begin{bmastersthesis}[author]
\bauthor{\bsnm{Donoho},~\bfnm{D.}\binits{D.}}
(\byear{1982}).
\btitle{Breakdown properties of multivariate location estimators}
\btype{Ph.{D}. {Q}ualifying {P}aper},
\bschool{Harvard University}.
\end{bmastersthesis}
\endbibitem

\bibitem[\protect\citeauthoryear{Donoho and Gasko}{1992}]{DonohoG92}
\begin{barticle}[author]
\bauthor{\bsnm{Donoho},~\bfnm{D.~L.}\binits{D.~L.}} \AND
  \bauthor{\bsnm{Gasko},~\bfnm{M.}\binits{M.}}
(\byear{1992}).
\btitle{Breakdown properties of location estimators based on halfspace depth
  and projected outlyingness}.
\bjournal{Annals of Statistics}
\bvolume{20}
\bpages{1803--1827}.
\end{barticle}
\endbibitem

\bibitem[\protect\citeauthoryear{D{\"{u}}mbgen}{1992}]{Duembgen92}
\begin{barticle}[author]
\bauthor{\bsnm{D{\"{u}}mbgen},~\bfnm{Lutz}\binits{L.}}
(\byear{1992}).
\btitle{Limit theorems for the simplicial depth}.
\bjournal{Statistics \& Probability Letters}
\bvolume{14}
\bpages{119--128}.
\end{barticle}
\endbibitem

\bibitem[\protect\citeauthoryear{Dyckerhoff}{2002}]{Dyckerhoff02}
\begin{bincollection}[author]
\bauthor{\bsnm{Dyckerhoff},~\bfnm{Rainer}\binits{R.}}
(\byear{2002}).
\btitle{Inference based on data depth}.
In \bbooktitle{Multivariate Dispersion, Central Regions and Depth: The Lift
  Zonoid Approach}
(\beditor{\bfnm{Karl}\binits{K.}~\bsnm{Mosler}}, ed.)
\bchapter{5},
\bpages{133--163}.
\bpublisher{Springer}, \baddress{New York}.
\end{bincollection}
\endbibitem

\bibitem[\protect\citeauthoryear{Dyckerhoff}{2004}]{Dyckerhoff04}
\begin{barticle}[author]
\bauthor{\bsnm{Dyckerhoff},~\bfnm{Rainer}\binits{R.}}
(\byear{2004}).
\btitle{Data depths satisfying the projection property}.
\bjournal{Allgemeines Statistisches Archiv}
\bvolume{88}
\bpages{163--190}.
\end{barticle}
\endbibitem

\bibitem[\protect\citeauthoryear{Dyckerhoff}{2017}]{Dyckerhoff17}
\begin{bunpublished}[author]
\bauthor{\bsnm{Dyckerhoff},~\bfnm{Rainer}\binits{R.}}
(\byear{2017}).
\btitle{Supplement to ``Convergence of depths and depth-trimmed regions''}.
\bnote{Unpublished manuscript}.
\end{bunpublished}
\endbibitem

\bibitem[\protect\citeauthoryear{Dyckerhoff and Mosler}{2011}]{DyckerhoffM11}
\begin{barticle}[author]
\bauthor{\bsnm{Dyckerhoff},~\bfnm{Rainer}\binits{R.}} \AND
  \bauthor{\bsnm{Mosler},~\bfnm{Karl}\binits{K.}}
(\byear{2011}).
\btitle{Weighted-mean trimming of multivariate data}.
\bjournal{Journal of Multivariate Analysis}
\bvolume{102}
\bpages{405--421}.
\end{barticle}
\endbibitem

\bibitem[\protect\citeauthoryear{Dyckerhoff and Mosler}{2012}]{DyckerhoffM12}
\begin{barticle}[author]
\bauthor{\bsnm{Dyckerhoff},~\bfnm{Rainer}\binits{R.}} \AND
  \bauthor{\bsnm{Mosler},~\bfnm{Karl}\binits{K.}}
(\byear{2012}).
\btitle{Weighted-mean regions of a probability distribution}.
\bjournal{Statistics \& Probability Letters}
\bvolume{82}
\bpages{318--325}.
\end{barticle}
\endbibitem

\bibitem[\protect\citeauthoryear{Eddy}{1985}]{Eddy85}
\begin{bincollection}[author]
\bauthor{\bsnm{Eddy},~\bfnm{W.~F.}\binits{W.~F.}}
(\byear{1985}).
\btitle{Ordering of multivariate data}.
In \bbooktitle{Compute Science and Statistics: Proceedings of the 16th
  Symposium on the Interface}
(\beditor{\bfnm{L.}\binits{L.}~\bsnm{Billard}}, ed.)
\bpages{25--30}.
\bpublisher{North-Holland}, \baddress{Amsterdam}.
\end{bincollection}
\endbibitem

\bibitem[\protect\citeauthoryear{He and Wang}{1997}]{HeW97}
\begin{barticle}[author]
\bauthor{\bsnm{He},~\bfnm{X.}\binits{X.}} \AND
  \bauthor{\bsnm{Wang},~\bfnm{G.}\binits{G.}}
(\byear{1997}).
\btitle{Convergence of depth contours for multivariate datasets}.
\bjournal{Annals of Statistics}
\bvolume{25}
\bpages{495--504}.
\end{barticle}
\endbibitem

\bibitem[\protect\citeauthoryear{Hoberg}{2000}]{Hoberg00}
\begin{binproceedings}[author]
\bauthor{\bsnm{Hoberg},~\bfnm{Richard}\binits{R.}}
(\byear{2000}).
\btitle{Cluster analysis based on data depth}.
In \bbooktitle{Data Analysis, Classification and Related Methods}
(\beditor{\bfnm{H.}\binits{H.}~\bsnm{Kiers}},
  \beditor{\bfnm{J.~P.}\binits{J.~P.}~\bsnm{Rasson}},
  \beditor{\bfnm{P.}\binits{P.}~\bsnm{Groenen}} \AND
  \beditor{\bfnm{M.}\binits{M.}~\bsnm{Schader}}, eds.)
\bpages{17--22}.
\bpublisher{Springer}, \baddress{Berlin}.
\end{binproceedings}
\endbibitem

\bibitem[\protect\citeauthoryear{Hoberg}{2003}]{Hoberg03}
\begin{bbook}[author]
\bauthor{\bsnm{Hoberg},~\bfnm{Richard}\binits{R.}}
(\byear{2003}).
\btitle{Clusteranalyse, Klassifikation und Datentiefe}.
\bpublisher{Eul}, \baddress{Lohmar}.
\end{bbook}
\endbibitem

\bibitem[\protect\citeauthoryear{Kim}{2000}]{Kim00}
\begin{barticle}[author]
\bauthor{\bsnm{Kim},~\bfnm{Jeankyung}\binits{J.}}
(\byear{2000}).
\btitle{Rate of convergence of depth contours: with application to a
  multivariate metrically trimmed mean}.
\bjournal{Statistics \& Probability Letters}
\bvolume{49}
\bpages{393--400}.
\end{barticle}
\endbibitem

\bibitem[\protect\citeauthoryear{Klein and Thompson}{1984}]{KleinT84}
\begin{bbook}[author]
\bauthor{\bsnm{Klein},~\bfnm{Erwin}\binits{E.}} \AND
  \bauthor{\bsnm{Thompson},~\bfnm{Anthony~C.}\binits{A.~C.}}
(\byear{1984}).
\btitle{Theory of Correspondences}.
\bpublisher{Wiley}, \baddress{New York}.
\end{bbook}
\endbibitem

\bibitem[\protect\citeauthoryear{Koshevoy and Mosler}{1997}]{KoshevoyM97}
\begin{barticle}[author]
\bauthor{\bsnm{Koshevoy},~\bfnm{Gleb}\binits{G.}} \AND
  \bauthor{\bsnm{Mosler},~\bfnm{Karl}\binits{K.}}
(\byear{1997}).
\btitle{Zonoid trimming for multivariate distributions}.
\bjournal{Annals of Statistics}
\bvolume{25}
\bpages{1998--2017}.
\end{barticle}
\endbibitem

\bibitem[\protect\citeauthoryear{Lange, Mosler and
  Mozharovskyi}{2014}]{LangeMM14}
\begin{barticle}[author]
\bauthor{\bsnm{Lange},~\bfnm{Tatjana}\binits{T.}},
  \bauthor{\bsnm{Mosler},~\bfnm{Karl}\binits{K.}} \AND
  \bauthor{\bsnm{Mozharovskyi},~\bfnm{Pavlo}\binits{P.}}
(\byear{2014}).
\btitle{Fast nonparametric classification based on data depth}.
\bjournal{Statistical Papers}
\bvolume{55}
\bpages{49--69}.
\end{barticle}
\endbibitem

\bibitem[\protect\citeauthoryear{Liu}{1988}]{Liu88}
\begin{barticle}[author]
\bauthor{\bsnm{Liu},~\bfnm{Regina~Y.}\binits{R.~Y.}}
(\byear{1988}).
\btitle{On a notion of simplicial depth}.
\bjournal{Proceedings of the National Academy of Sciences of the USA}
\bvolume{85}
\bpages{1732--1734}.
\end{barticle}
\endbibitem

\bibitem[\protect\citeauthoryear{Liu}{1990}]{Liu90}
\begin{barticle}[author]
\bauthor{\bsnm{Liu},~\bfnm{R.~Y.}\binits{R.~Y.}}
(\byear{1990}).
\btitle{On a notion of data depth based on random simplices}.
\bjournal{Annals of Statistics}
\bvolume{18}
\bpages{405--414}.
\end{barticle}
\endbibitem

\bibitem[\protect\citeauthoryear{Liu}{1992}]{Liu92}
\begin{bincollection}[author]
\bauthor{\bsnm{Liu},~\bfnm{Regina~Y.}\binits{R.~Y.}}
(\byear{1992}).
\btitle{Data depth and multivariate rank tests}.
In \bbooktitle{$L_1$-Statistical Analysis and Related Methods}
(\beditor{\bfnm{Y.}\binits{Y.}~\bsnm{Dodge}}, ed.)
\bpublisher{North-Holland}, \baddress{Amsterdam}.
\end{bincollection}
\endbibitem

\bibitem[\protect\citeauthoryear{Liu}{1995}]{Liu95}
\begin{barticle}[author]
\bauthor{\bsnm{Liu},~\bfnm{Regina~Y.}\binits{R.~Y.}}
(\byear{1995}).
\btitle{Control charts for multivariate processes}.
\bjournal{Journal of the American Statistical Association}
\bvolume{90}
\bpages{1380--1387}.
\end{barticle}
\endbibitem

\bibitem[\protect\citeauthoryear{Liu, Parelius and Singh}{1999}]{LiuPS99}
\begin{barticle}[author]
\bauthor{\bsnm{Liu},~\bfnm{R.~Y.}\binits{R.~Y.}},
  \bauthor{\bsnm{Parelius},~\bfnm{J.~M.}\binits{J.~M.}} \AND
  \bauthor{\bsnm{Singh},~\bfnm{K.}\binits{K.}}
(\byear{1999}).
\btitle{Multivariate analysis by data depth: Descriptive statistics, graphics
  and inference}.
\bjournal{Annals of Statistics}
\bvolume{27}
\bpages{783--858}.
\end{barticle}
\endbibitem

\bibitem[\protect\citeauthoryear{Liu and Singh}{1993}]{LiuS93}
\begin{barticle}[author]
\bauthor{\bsnm{Liu},~\bfnm{Regina~Y.}\binits{R.~Y.}} \AND
  \bauthor{\bsnm{Singh},~\bfnm{Kesar}\binits{K.}}
(\byear{1993}).
\btitle{A quality index based on data depth and multivariate rank tests}.
\bjournal{Journal of the American Statistical Association}
\bvolume{88}
\bpages{252--260}.
\end{barticle}
\endbibitem

\bibitem[\protect\citeauthoryear{Mahalanobis}{1936}]{Mahalanobis36}
\begin{barticle}[author]
\bauthor{\bsnm{Mahalanobis},~\bfnm{P.~C.}\binits{P.~C.}}
(\byear{1936}).
\btitle{On the generalized distance in statistics}.
\bjournal{Procedings of the National Academy of India}
\bvolume{12}
\bpages{49--55}.
\end{barticle}
\endbibitem

\bibitem[\protect\citeauthoryear{Mass{\'{e}}}{2002}]{Masse02}
\begin{barticle}[author]
\bauthor{\bsnm{Mass{\'{e}}},~\bfnm{Jean-Claude}\binits{J.-C.}}
(\byear{2002}).
\btitle{Asymptotics for the Tukey median}.
\bjournal{Journal of Multivariate Analysis}
\bvolume{81}
\bpages{286--300}.
\end{barticle}
\endbibitem

\bibitem[\protect\citeauthoryear{Mass{\'{e}}}{2004}]{Masse04}
\begin{barticle}[author]
\bauthor{\bsnm{Mass{\'{e}}},~\bfnm{Jean-Claude}\binits{J.-C.}}
(\byear{2004}).
\btitle{Asymptotics for the Tukey depth process, with an application to a
  multivariate trimmed mean}.
\bjournal{Bernoulli}
\bvolume{10}
\bpages{397--419}.
\end{barticle}
\endbibitem

\bibitem[\protect\citeauthoryear{Mass{\'{e}} and Theodorescu}{1994}]{MasseT94}
\begin{barticle}[author]
\bauthor{\bsnm{Mass{\'{e}}},~\bfnm{J.~C.}\binits{J.~C.}} \AND
  \bauthor{\bsnm{Theodorescu},~\bfnm{R.}\binits{R.}}
(\byear{1994}).
\btitle{Halfplane trimming for bivariate distributions}.
\bjournal{The American Statistician}
\bvolume{48}
\bpages{188--202}.
\end{barticle}
\endbibitem

\bibitem[\protect\citeauthoryear{Mosler}{2002}]{Mosler02}
\begin{bbook}[author]
\bauthor{\bsnm{Mosler},~\bfnm{Karl}\binits{K.}}
(\byear{2002}).
\btitle{Multivariate Dispersion, Central Regions and Depth}.
\bpublisher{Springer}, \baddress{New York}.
\end{bbook}
\endbibitem

\bibitem[\protect\citeauthoryear{Mosler and Bazovkin}{2014}]{MoslerB14}
\begin{barticle}[author]
\bauthor{\bsnm{Mosler},~\bfnm{Karl}\binits{K.}} \AND
  \bauthor{\bsnm{Bazovkin},~\bfnm{Pavel}\binits{P.}}
(\byear{2014}).
\btitle{Stochastic linear programming with a distortion risk constraint}.
\bjournal{OR Spectrum}
\bvolume{36}
\bpages{949--969}.
\end{barticle}
\endbibitem

\bibitem[\protect\citeauthoryear{Mosler and Hoberg}{2006}]{MoslerH06}
\begin{bincollection}[author]
\bauthor{\bsnm{Mosler},~\bfnm{Karl}\binits{K.}} \AND
  \bauthor{\bsnm{Hoberg},~\bfnm{Richard}\binits{R.}}
(\byear{2006}).
\btitle{Data analysis and classification with the zonoid depth}.
In \bbooktitle{Data Depth: Robust Multivariate Analysis, Computational Geometry
  and Applications}
(\beditor{\bfnm{R.}\binits{R.}~\bsnm{Liu}},
  \beditor{\bfnm{R.}\binits{R.}~\bsnm{Serfling}} \AND
  \beditor{\bfnm{D.}\binits{D.}~\bsnm{Souvaine}}, eds.)
\bpages{49--59}.
\bpublisher{American Mathematical Society}, \baddress{Providence RI}.
\end{bincollection}
\endbibitem

\bibitem[\protect\citeauthoryear{Nolan}{1992}]{Nolan92}
\begin{barticle}[author]
\bauthor{\bsnm{Nolan},~\bfnm{D.}\binits{D.}}
(\byear{1992}).
\btitle{Asymptotics for multivariate trimming}.
\bjournal{Stochastic Processes and Their Applications}
\bvolume{42}
\bpages{157--169}.
\end{barticle}
\endbibitem

\bibitem[\protect\citeauthoryear{Rousseeuw and Ruts}{1999}]{RousseeuwR99}
\begin{barticle}[author]
\bauthor{\bsnm{Rousseeuw},~\bfnm{P.~J.}\binits{P.~J.}} \AND
  \bauthor{\bsnm{Ruts},~\bfnm{Ida}\binits{I.}}
(\byear{1999}).
\btitle{The depth function of a population distribution}.
\bjournal{Metrika}
\bvolume{49}
\bpages{213--244}.
\end{barticle}
\endbibitem

\bibitem[\protect\citeauthoryear{Serfling}{2002a}]{Serfling02b}
\begin{barticle}[author]
\bauthor{\bsnm{Serfling},~\bfnm{Robert}\binits{R.}}
(\byear{2002}a).
\btitle{Generalized quantile processes based on multivariate depth functions,
  with applications in nonparametric multivariate analysis}.
\bjournal{Journal of Multivariate Analysis}
\bvolume{83}
\bpages{232--247}.
\end{barticle}
\endbibitem

\bibitem[\protect\citeauthoryear{Serfling}{2002b}]{Serfling02a}
\begin{barticle}[author]
\bauthor{\bsnm{Serfling},~\bfnm{Robert}\binits{R.}}
(\byear{2002}b).
\btitle{Quantile functions for multivariate analysis: approaches and
  applications}.
\bjournal{Statistica Neerlandica}
\bvolume{56}
\bpages{214--232}.
\end{barticle}
\endbibitem

\bibitem[\protect\citeauthoryear{Singh}{1991}]{Singh91}
\begin{btechreport}[author]
\bauthor{\bsnm{Singh},~\bfnm{K.}\binits{K.}}
(\byear{1991}).
\btitle{A notion of majority depth}
\btype{Technical Report},
\binstitution{Rutgers University, Department of Statistics}.
\end{btechreport}
\endbibitem

\bibitem[\protect\citeauthoryear{Stahel}{1981}]{Stahel81}
\begin{bphdthesis}[author]
\bauthor{\bsnm{Stahel},~\bfnm{W.~A.}\binits{W.~A.}}
(\byear{1981}).
\btitle{Robuste Sch{\"{a}}tzungen: Infinitesimale Optimalit{\"{a}}t und
  Sch{\"{a}}tzungen von Kovarianzmatrizen}
\btype{PhD thesis},
\bschool{ETH Z{\"{u}}rich}.
\end{bphdthesis}
\endbibitem

\bibitem[\protect\citeauthoryear{Tukey}{1975}]{Tukey75}
\begin{binproceedings}[author]
\bauthor{\bsnm{Tukey},~\bfnm{J.~W.}\binits{J.~W.}}
(\byear{1975}).
\btitle{Mathematics and picturing data}.
In \bbooktitle{Proceedings of the 1974 International Congress of
  Mathematicians, Vancouver}
(\beditor{\bfnm{R.}\binits{R.}~\bsnm{James}}, ed.)
\bvolume{2}
\bpages{523--531}.
\end{binproceedings}
\endbibitem

\bibitem[\protect\citeauthoryear{Yeh and Singh}{1997}]{YehS97}
\begin{barticle}[author]
\bauthor{\bsnm{Yeh},~\bfnm{A.~B.}\binits{A.~B.}} \AND
  \bauthor{\bsnm{Singh},~\bfnm{Kesar}\binits{K.}}
(\byear{1997}).
\btitle{Balanced confidence regions based on Tukey`s depth and the bootstrap}.
\bjournal{Journal of the Royal Statistical Society, Series B}
\bvolume{59}
\bpages{639--652}.
\end{barticle}
\endbibitem

\bibitem[\protect\citeauthoryear{Zuo}{2003}]{Zuo03}
\begin{barticle}[author]
\bauthor{\bsnm{Zuo},~\bfnm{Yijun}\binits{Y.}}
(\byear{2003}).
\btitle{Projection based depth functions and associated medians}.
\bjournal{Annals of Statistics}
\bvolume{31}
\bpages{1460--1490}.
\end{barticle}
\endbibitem

\bibitem[\protect\citeauthoryear{Zuo}{2006}]{Zuo06}
\begin{barticle}[author]
\bauthor{\bsnm{Zuo},~\bfnm{Yijun}\binits{Y.}}
(\byear{2006}).
\btitle{Multidimensional trimming based on projection depth}.
\bjournal{Annals of Statistics}
\bvolume{34}
\bpages{2211-2251}.
\end{barticle}
\endbibitem

\bibitem[\protect\citeauthoryear{Zuo and Serfling}{2000a}]{ZuoS00a}
\begin{barticle}[author]
\bauthor{\bsnm{Zuo},~\bfnm{Yijun}\binits{Y.}} \AND
  \bauthor{\bsnm{Serfling},~\bfnm{Robert}\binits{R.}}
(\byear{2000}a).
\btitle{General notions of statistical depth function}.
\bjournal{Annals of Statistics}
\bvolume{28}
\bpages{461--482}.
\end{barticle}
\endbibitem

\bibitem[\protect\citeauthoryear{Zuo and Serfling}{2000b}]{ZuoS00b}
\begin{barticle}[author]
\bauthor{\bsnm{Zuo},~\bfnm{Yijun}\binits{Y.}} \AND
  \bauthor{\bsnm{Serfling},~\bfnm{Robert}\binits{R.}}
(\byear{2000}b).
\btitle{Structural properties and convergence results for contours of sample
  statistical depth functions}.
\bjournal{Annals of Statistics}
\bvolume{28}
\bpages{483--499}.
\end{barticle}
\endbibitem

\end{thebibliography}
\end{document}


\begin{frontmatter}
\title{Supplement to ``Convergence of depths and depth-trimmed regions''}
\runtitle{Supplement to ''Convergence of depths''}
\begin{aug}
\author{\fnms{Rainer} \snm{Dyckerhoff}\ead[label=e1]{rainer.dyckerhoff@statistik.uni-koeln.de}
\ead[label=u1,url]{http://www.wisostat.uni-koeln.de/institut/professoren/dyckerhoff/}}
\runauthor{R. Dyckerhoff}
\affiliation{University of Cologne\thanksmark{m1}}
\address{Institute of Econometrics and Statistics\\
University of Cologne\\
Albertus-Magnus-Platz\\
50923 Cologne\\
Germany\\
\printead{e1}\\
\printead{u1}}
\end{aug}
\end{frontmatter}
\manuallabel{th.PktTPktZ}{4.2}
\manuallabel{th.PktZPktT}{4.3}
\manuallabel{th.KomTKomZ}{4.5}
\manuallabel{th.KomZGlmT}{4.6}
\manuallabel{th.PktZKomZ}{4.7}
\manuallabel{th.PktTGlmT}{4.8}
\section{Introduction}
In this supplement we present some examples to show that
\begin{enumerate}
\item without the assumption `strictly monotone for $P$' 
Theorems~\ref{th.PktTPktZ}, \ref{th.KomTKomZ}, and \ref{th.PktZKomZ}
are in general false,
\item without the assumption `continuous for $P$' 
Theorems~\ref{th.PktZPktT},  \ref{th.KomZGlmT}, and  \ref{th.PktTGlmT}
are in general false.
\end{enumerate}


\section{Examples}
\begin{example}['Strict monotonicity' is needed in Theorems~\ref{th.PktTPktZ} and \ref{th.KomTKomZ}]\label{ex.1}
Let $Q_1=U\bigl([-3,-2]\union[2,3]\bigr)$ be the uniform distribution
on the union of the two intervals $[-3,-2]$ and $[3,2]$. Let further
$Q_2=U\bigl([-1,1]\bigr)$. Now, for $n\in\IN$ let
\[
P_n=\frac12\left(1+\frac{(-1)^n}{n}\right)Q_1+\frac12\left(1-\frac{(-1)^n}{n}\right)Q_2
\]
and finally
\[
P_0=\frac12Q_1+\frac12Q_2\,.
\]
The following figure shows the functions $x\mapsto P_0\bigl((-\infty,x]\bigr)$, 
$x\mapsto P_0\bigl([x,\infty)\bigr)$ as well as the halfspace depth (in blue), 
\[
x\mapsto\HD(x\,|\,P_0)=\min\{P_0\bigl((-\infty,x]), P_0([x,\infty)\bigr)\}\,.
\]
\[
\psset{xunit=1.3cm,yunit=3cm}\footnotesize
\begin{pspicture}(-3.5,-0.2)(3.5,1.2)
\psaxes[arrows=->,Dx=1,Dy=0.25,labels=all](0,0)(-3.5,-0.2)(3.5,1.2)
\psline(-3.5,0)(-3,0)(-2,0.25)(-1,0.25)(1,0.75)(2,0.75)(3,1)(3.5,1)
\psline(-3.5,1)(-3,1)(-2,0.75)(-1,0.75)(1,0.25)(2,0.25)(3,0)(3.5,0)
\psline[linewidth=1.5pt,linecolor=blue](-3.5,0)(-3,0)(-2,0.25)(-1,0.25)(0,0.5)(1,0.25)(2,0.25)(3,0)(3.5,0)
\uput[u](1.5,0.75){$P_0\bigl((-\infty,x]\bigr)$}
\uput[u](-1.5,0.75){$P_0\bigl([x,\infty)\bigr)$}
\uput[u](1.5,0.25){$\HD(x\,|\,P_0)$}
\uput[r](3.5,0){$x$}
\end{pspicture}
\]
Note, that the halfspace depth is not strictly monotone for $P_0$,
since
\[
\HD_\frac14(P_0)=[-2,2]\ne [-1,1]=\cl\left\{z\in\IR\,|\,\HD(z\,|\,P_0)>\frac14\right\}\,.
\]
In the following figure the halfspace depth w.r.t. $P_0$ (blue line) as well as 
the halfspace depths w.r.t. $P_9$ (red line below blue line) and $P_{10}$ 
(red line above blue line) are drawn. 
\[
\psset{xunit=1.3cm,yunit=6cm}\footnotesize
\begin{pspicture}(-3.5,-0.05)(3.5,0.55)
\psaxes[arrows=->,Dx=1,Dy=0.25,labels=all](0,0)(-3.5,-0.05)(3.5,0.55)
\psline[linecolor=red](-3.5,0)(-3,0)(-2,0.275)(-1,0.275)(0,0.5)(1,0.275)(2,0.275)(3,0)(3.5,0)
\psline[linecolor=blue](-3.5,0)(-3,0)(-2,0.25)(-1,0.25)(0,0.5)(1,0.25)(2,0.25)(3,0)(3.5,0)
\psline[linecolor=red](-3.5,0)(-3,0)(-2,0.2222)(-1,0.2222)(0,0.5)(1,0.2222)(2,0.2222)(3,0)(3.5,0)
\uput[u](1.7,0.26){$\HD(x\,|\,P_{10})$}
\uput[d](1.3,0.24){$\HD(x\,|\,P_{9})$}
\uput[r](3.5,0){$x$}
\end{pspicture}
\]   
It is easy to see that
\begin{align*}
\HD(z\,|\,P_n) \le \HD(z\,|\,P_0)\ &\text{for all $z$, if $n$ is odd,}\\
\HD(z\,|\,P_n) \ge \HD(z\,|\,P_0)\ &\text{for all $z$, if $n$ is even.}
\end{align*}
Further,
\[
\sup_{z\in\IR}\bigl|\HD(z\,|\,P_n)-\HD(z\,|\,P_0)\bigr|=\frac1{4n}\,.
\]
Therefore,
\[
\lim_{n\to\infty}\sup_{z\in\IR}\bigl|\HD(z\,|\,P_n)-\HD(z\,|\,P_0)\bigr|=0\,,
\]
i.e., \UniD and a-fortiori \ComD and \PtwD hold.

On the other hand we have
\begin{align*}
\HD(1\,|\,P_n)=\HD(2\,|\,P_n)&< \frac14 = \HD(2\,|\,P_0)\quad\text{if $n$ is odd,}\\
\HD(1\,|\,P_n)=\HD(2\,|\,P_n)&> \frac14 = \HD(2\,|\,P_0)\quad\text{if $n$ is even.}
\end{align*} 
Therefore and because of symmetry of the depth functions, for $m$ odd and $n$ 
even holds  
\[
\HD_{0.25}(P_m)\subset(-1,1)
\subset\HD_{0.25}(P_0)=[-2,2]
\subset\HD_{0.25}(P_n)\,.
\]
This shows that although the depth functions converge uniformly, the trimmed 
regions $\HD_{0.25}(P_n)$ do not converge, neither in the Hausdorff sense nor
in the set-theoretic sense. Thus, neither \ComR nor \PtwR hold.
In other words, without the assumption of 'strict monotonicity for $P$',
Theorems~\ref{th.PktTPktZ} and \ref{th.KomTKomZ} do in general not hold. 
\end{example}


\begin{example}['Strict monotonicity' is needed in Theorem~\ref{th.PktZKomZ}]\label{ex.2}
This example is similar to the previous one with the main difference that
the alternating term $(-1)^n$ is missing. So, let 
$Q_1=U\bigl([-3,-2]\union[2,3]\bigr)$ and $Q_2=U\bigl([-1,1]\bigr)$ as in
Example~\ref{ex.1}. Now, for $n\in\IN$ let
\[
P_n=\frac12\left(1+\frac{1}{n}\right)Q_1+\frac12\left(1-\frac{1}{n}\right)Q_2
\]
and finally $P_0=\frac12Q_1+\frac12Q_2$ as in Example~\ref{ex.1}.
In the following figure the halfspace depth w.r.t. $P_0$ (blue line) as well as 
the halfspace depth w.r.t. $P_{10}$ (red line above blue line) are shown. 
\[
\psset{xunit=1.3cm,yunit=6cm}\footnotesize
\begin{pspicture}(-3.5,-0.05)(3.5,0.55)
\psaxes[arrows=->,Dx=1,Dy=0.25,labels=all](0,0)(-3.5,-0.05)(3.5,0.55)
\psline[linecolor=red](-3.5,0)(-3,0)(-2,0.275)(-1,0.275)(0,0.5)(1,0.275)(2,0.275)(3,0)(3.5,0)
\psline[linecolor=blue](-3.5,0)(-3,0)(-2,0.25)(-1,0.25)(0,0.5)(1,0.25)(2,0.25)(3,0)(3.5,0)
\uput[u](1.7,0.26){$\HD(x\,|\,P_{10})$}
\uput[d](1.3,0.25){$\HD(x\,|\,P_0)$}
\uput[r](3.5,0){$x$}
\end{pspicture}
\]   
As was shown in Example~\ref{ex.1}, the halfspace depth is not strictly monotone 
for $P_0$.

Since in this example the depth function w.r.t. $P_n$ lies always above the
depth function w.r.t. $P_0$, the oscillating behavior of the trimmed regions
$\HD_{0.25}(P_n)$ in the previous example does not occur. Instead, for each
$\alpha$ (even for $\alpha=0.25$) the trimmed regions $\HD_\alpha(P_n)$ converge
to $\HD_\alpha(P_0)$ in the Hausdorff metric. Therefore \PtwR holds.

However, \ComR does not hold. This can be seen as follows. For $n$ arbitrarily,
choose $\alpha_n$ such that
\[
\frac14<\alpha_n<\frac14\left(1+\frac1n\right)\,.
\]
Then,
\[
\HD(1\,|\,P_0)=\frac14<\alpha_n<\frac14\left(1+\frac1n\right)=\HD(2\,|\,P_n)\,.
\]
Therefore and because of symmetry of the depth function,
\[
\HD_{\alpha_n}(P_0)\subset(-1,1)\subset[-2,2]\subset \HD_{\alpha_n}(P_n)\,.
\]
This shows that $\delta_H\bigl(\HD_{\alpha_n}(P_n),\HD_{\alpha_n}(P_0)\bigr)\ge 1$.
Since this holds for every $n$, it follows
\[
\sup_{\alpha\in [0.1,0.5]}\delta_H\bigl(\D_\alpha(P_n),\D_\alpha(P_0)\bigr)
\ge 1\quad\text{for all $n$,}
\]
i.e., there is a compact interval $A\subset\bigl(0,\alpha_{\max}(P_0)\bigr)$
on which the trimmed regions do not converge uniformly.  
In other words, without the assumption of 'strict monotonicity for $P$',
\PtwR does not imply \ComR and thus, Theorem~\ref{th.PktZKomZ} does in general 
not hold. 
\end{example}


\begin{example}['Continuity' is needed in Theorems~\ref{th.PktZPktT} and \ref{th.KomZGlmT}]\label{ex.3}
For $n\in\IN_0$ define
\[
a_n=\begin{cases}
1+\frac{(-1)^n}{n+1}\,,&\text{if $n\in\IN$,}\\
1\,&\text{if $n=0$.}
\end{cases}
\]
Now, let $Q_{1,n}=U\bigl([-2,-a_n]\union[a_n,2]\bigr)$ be the uniform 
distribution on the union of the intervals $[-2,-a_n]$ and $[a_n,2]$,
$Q_{2,n}=U\bigl([-a_n,a_n]\bigr)$ and 
$Q_{3,n}=0.5\epsilon_{-a_n}+0.5\epsilon_{a_n}$,
where $\epsilon_x$ denotes the one-point measure on $x$.
For $n\in\IN_0$ consider the probability measures
\[
P_n=0.3\,Q_{1,n}+0.3\,Q_{2,n}+0.4\,Q_{3,n}\,.
\]
The following figure shows the functions $x\mapsto P_0\bigl((-\infty,x]\bigr)$, 
$x\mapsto P_0\bigl([x,\infty)\bigr)$ as well as the halfspace depth (in blue), 
\[
x\mapsto\HD(x\,|\,P_0)=\min\{P_0\bigl((-\infty,x]),P_0([x,\infty)\bigr)\}\,.
\]
\[
\psset{xunit=2cm,yunit=3cm}\footnotesize
\begin{pspicture}(-2.5,-0.2)(2.5,1.2)
\psaxes[arrows=->,Dx=1,labels=all](0,0)(-2.5,-0.2)(2.5,1.2)
\uput[l](-4pt,0.15){$0.15$}\psline[linewidth=0.6pt](-4pt,0.15)(4pt,0.15)
\uput[l](-4pt,0.35){$0.35$}\psline[linewidth=0.6pt](-4pt,0.35)(4pt,0.35)
\uput[l](-4pt,0.65){$0.65$}\psline[linewidth=0.6pt](-4pt,0.65)(4pt,0.65)
\uput[l](-4pt,0.85){$0.85$}\psline[linewidth=0.6pt](-4pt,0.85)(4pt,0.85)
\psline(-2.5,0)(-2,0)(-1,0.15)\psline[arrows=*-](-1,0.35)(1,0.65)\psline[arrows=*-](1,0.85)(2,1)(2.5,1)
\psline[arrows=-*](-2.5,1)(-2,1)(-1,0.85)\psline[arrows=-*](-1,0.65)(1,0.35)\psline(1,0.15)(2,0)(2.5,0)
\psset{linewidth=1.5pt,linecolor=blue}
\psline(-2.5,0)(-2,0)(-1,0.15)\psline[arrows=*-*](-1,0.35)(0,0.5)(1,0.35)\psline(1,0.15)(2,0)(2.5,0)
\uput[dl](2.5,0.95){$P_0\bigl((-\infty,x]\bigr)$}
\uput[dr](-2.5,0.95){$P_0\bigl([x,\infty)\bigr)$}
\uput[ul](2.5,0.05){$\HD(x\,|\,P_0)$}
\uput[r](2.5,0){$x$}
\end{pspicture}
\]
Obviously, the halfspace depth is not continuous for $P_0$. 

In the following figure the halfspace depth w.r.t. $P_0$ (blue line) as well as 
the halfspace depths w.r.t. $P_8$ (red line above blue line) and $P_9$ 
(red line below blue line) are drawn. 
\[
\psset{xunit=2cm,yunit=6cm}\footnotesize
\begin{pspicture}(-2.5,-0.05)(2.5,0.55)
\psaxes[arrows=->,Dx=1,Dy=0.5,labels=all](0,0)(-2.5,-0.05)(2.5,0.55)
\uput[l](-4pt,0.15){$0.15$}\psline[linewidth=0.6pt](-4pt,0.15)(4pt,0.15)
\uput[l](-4pt,0.35){$0.35$}\psline[linewidth=0.6pt](-4pt,0.35)(4pt,0.35)
\psset{linecolor=red}
\psline(-2.5,0)(-2,0)(-1.1111,0.15)\psline[arrows=*-*](-1.1111,0.35)(0,0.5)(1.1111,0.35)\psline(1.1111,0.15)(2,0)(2.5,0)
\psset{linecolor=blue}
\psline(-2.5,0)(-2,0)(-1,0.15)\psline[arrows=*-*](-1,0.35)(0,0.5)(1,0.35)\psline(1,0.15)(2,0)(2.5,0)
\psset{linecolor=red}
\psline(-2.5,0)(-2,0)(-0.9,0.15)\psline[arrows=*-*](-0.9,0.35)(0,0.5)(0.9,0.35)\psline(0.9,0.15)(2,0)(2.5,0)
\uput[ur](1.4,0.08){$\HD(x\,|\,P_{8})$}
\uput[dl](1.1,0.12){$\HD(x\,|\,P_{9})$}
\uput[r](2.5,0){$x$}
\end{pspicture}
\]   
It is easy to see that
\begin{align*}
\HD(z\,|\,P_n) \le \HD(x\,|\,P_0)\ &\text{for all $z$, if $n$ is odd,}\\
\HD(z\,|\,P_n) \ge \HD(x\,|\,P_0)\ &\text{for all $z$, if $n$ is even.}
\end{align*}
Further,
\[
\sup_{\alpha\in [0,1/2]}\delta_H\bigl(\D_\alpha(P_n),\D_\alpha(P_0)\bigr)=\frac1{n+1}\,.
\]
Therefore,
\[
\lim_{n\to\infty}\sup_{\alpha\in [0,1/2]}\delta_H\bigl(\D_\alpha(P_n),\D_\alpha(P_0)\bigr)=0\,,
\]
i.e., \ComR and a-fortiori \PtwR hold.

On the other hand we have
\[
\HD(1\,|\,P_m)<0.15 = \lim_{x\searrow1}\HD(x\,|\,P_0)
< \HD(1\,|\,P_0)=0.35 < \HD(1\,|\,P_n)
\]
for $m$ odd and $n$ even. This shows that although the trimmed regions
converge uniformly, the depth functions $\HD(x\,|\,P_n)$ do not converge
for $x=1$. Thus, neither \UniD nor \PtwD hold.
In other words, without the assumption of 'continuity for $P$',
Theorems~\ref{th.PktZPktT} and \ref{th.KomZGlmT} do in general not hold. 
\end{example}


\begin{example}['Continuity' is needed in Theorem~\ref{th.PktTGlmT}]\label{ex.4}
This example is similar to the previous one with the main difference that
the alternating term $(-1)^n$ is missing.
So, for $n\in\IN_0$ let
\[
a_n=\begin{cases}
1+\frac{1}{n+1}\,,&\text{if $n\in\IN$,}\\
1\,&\text{if $n=0$.}
\end{cases}
\]
Define $Q_{1,n}=U\bigl([-2,-a_n]\union[a_n,2]\bigr)$, 
$Q_{2,n}=U\bigl([-a_n,a_n]\bigr)$ and 
$Q_{3,n}=0.5\epsilon_{-a_n}+0.5\epsilon_{a_n}$ as in Example~\ref{ex.3}.
Finally, for $n\in\IN_0$ consider the probability measures
\[
P_n=0.3\,Q_{1,n}+0.3\,Q_{2,n}+0.4\,Q_{3,n}\,.
\]
As in Example~\ref{ex.3}, the halfspace depth is not continuous for $P_0$. 

In the following figure the halfspace depth w.r.t. $P_0$ (blue line) as well as 
the halfspace depth w.r.t. $P_8$ (red line above blue line) are shown. 
\[
\psset{xunit=2cm,yunit=6cm}\footnotesize
\begin{pspicture}(-2.5,-0.05)(2.5,0.55)
\psaxes[arrows=->,Dx=1,Dy=0.55,labels=all](0,0)(-2.5,-0.05)(2.5,0.55)
\uput[l](-4pt,0.15){$0.15$}\psline[linewidth=0.6pt](-4pt,0.15)(4pt,0.15)
\uput[l](-4pt,0.35){$0.35$}\psline[linewidth=0.6pt](-4pt,0.35)(4pt,0.35)
\psset{linecolor=red}
\psline(-2.5,0)(-2,0)(-1.1111,0.15)\psline[arrows=*-*](-1.1111,0.35)(0,0.5)(1.1111,0.35)\psline(1.1111,0.15)(2,0)(2.5,0)
\psset{linecolor=blue}
\psline(-2.5,0)(-2,0)(-1,0.15)\psline[arrows=*-*](-1,0.35)(0,0.5)(1,0.35)\psline(1,0.15)(2,0)(2.5,0)
\uput[ur](1.4,0.08){$\HD(x\,|\,P_{8})$}
\uput[dl](1.1,0.12){$\HD(x\,|\,P_{0})$}
\uput[r](2.5,0){$x$}
\end{pspicture}
\]
Since in this example the depth function w.r.t. $P_n$ lies always above the
depth function w.r.t. $P$, the oscillating behavior of the depth $\HD(1\,|\,P_n)$ 
in the previous example does not occur. Instead, for each $x$ (even for $x=\pm1$) 
the depth $\HD(x\,|\,P_n)$ converges to $\HD(x\,|\,P_0)$. Therefore, \PtwD holds.

However, \UniD does not hold. This can be seen as follows. For $n$ arbitrarily,
choose $x_n$ such that $1<x_n<a_n$. Since $x_n$ lies between the jumps of
$\HD(x\,|\,P_0)$ and $\HD(x\,|\,P_n)$ it follows
\[
\HD(x_n\,|\,P_0)<0.15<0.35<\HD(x_n\,|\,P_n)\,.
\]
This shows that $\left|\HD(x_n\,|\,P_0)-\HD(x_n\,|\,P_n)\right|>0.2$.
Since this holds for every $n$, it holds
\[
\sup_{z\in\IR}\left|\HD(z\,|\,P_0)-\HD(z\,|\,P_n)\right|>0.2\quad\text{for all $n$,}
\]
i.e., the depth functions do not converge uniformly. In other words, without 
the assumption of 'continuity for $P$', \PtwD does not imply \UniD and thus, 
Theorem~\ref{th.PktTGlmT} does in general 
not hold. 
\end{example}